\documentclass[leqno]{article}

\def\ad{\mathop{\rm ad}}
\def\im{\mathop{\rm Im}}
\def\re{\mathop{\rm Re}}
\def\tr{\mathop{\rm tr}}

\begin{document}

\title{Notes on Lie Algebras and Lie Groups}

\author{Stephen William Semmes	\\
	Rice University		\\
	Houston, Texas}

\date{}

\maketitle

\tableofcontents

\section{Groups}
\setcounter{equation}{0}

	A group is a nonempty set $G$ with a distinguished element $e$
and a binary operation such that $e$ is the identity element for the
binary operation, the group operation satisfies the associative law,
and every element of the group has an inverse.  If also the group
operation satisfies the commutative law, then we say that $G$ is a
commutative or abelian group.

	For instance, if $E$ is a nonempty set, then the collection of
one-to-one mappings of $E$ onto itself forms a group, using the
identity mapping which fixes each element of $E$ as the identity
element of the group, and using composition of mappings as the group
operation.

	If $G$ is a group and $H$ is a subset of $G$ which contains
the identity element $e$, is closed under the group operation, and
contains the inverse of all of its elements, then $H$ is a subgroup of
$G$.  Basically this means that $H$ is a group too with the same
identity element and the restriction of the group operation to $H$.

	If $G_1$, $G_2$ are groups and $\phi$ is a mapping from $G_1$
to $G_2$, then $\phi$ is a group homomorphism if $\phi$ maps the
identity element of $G_1$ to the identity element of $G_2$ and if the
group operation on $G_1$ corresponds to the group operation on $G_2$
under $\phi$.  If $\phi$ is a one-to-one mapping of $G_1$ onto $G_2$,
then the inverse mapping $\phi^{-1} : G_2 \to G_1$ is automatically a
homomorphism too, and we say that $\phi$ is an isomorphism from $G_1$
onto $G_2$.

	If $G_1$, $G_2$ are groups and $\phi$ is a homomorphism from
$G_1$ into $G_2$, then the image of $\phi$, consisting of the points
$\phi(x)$ with $x \in G_1$, is a subgroup of $G_2$.  The kernel of
$\phi$ is the subgroup of $G_1$ consisting of elements which are
mapped by $\phi$ to the identity element of $G_2$.  One can check that
$\phi$ is one-to-one if and only if the kernel of $\phi$ is trivial,
which is to say that it contains only the identity element of $G_1$.

	A subgroup $H$ of a group $G$ is said to be normal if $g \, h
\, g^{-1} \in H$ for all $g \in G$ and $h \in H$.  The kernel of a
homomorphism from $G$ into another group is a normal subgroup, and
every subgroup of an abelian group is normal.

	Suppose that $G$ is a group and that $H$ is a subgroup.  We
can define equivalence relations $\sim_l$, $\sim_r$ on $G$ by saying
that $x \sim_l y$ if there is an $h \in H$ such that $y = x \, h$, and
that $x \sim_r y$ if there is an $h \in H$ such that $y = h \, x$.
These define equivalence relations on $G$, which is to say that they
are reflexive, symmetric, and transitive relations on $G$.

	These two equivalence relations are trivially the same if $G$
is an abelian group.  They are also the same if $H$ is a normal
subgroup of $G$.

	The equivalence classes in $G$ associated to the equivalence
relations $\sim_l$, $\sim_r$ are called the left and right cosets of
$H$ in $G$, respectively.  If $x \in G$, then $x \, H$ is defined to
be the subset of $G$ consisting of $x \, h$, $h \in G$, and $H \, x$
to be the subset of $G$ consisting of $h \, x$, $h \in H$.  These are
the left and right cosets of $H$ in $G$ containing $x$, which are the
same as the set of $y \in G$ such that $y \sim_l x$ and $y \sim_r x$,
respectively.

	The spaces of left and right cosets of $H$ in $G$ are denoted
$G / H$ and $H \backslash G$, respectively.  There are canonical
mappings from $G$ onto the coset spaces, which take an element of $G$
to the coset containing it.

	If $g \in G$, then we get a mapping from $G / H$ to itself by
sending a left coset $x \, H$ to $g \, x \, H$.  Similarly, it is
convenient to define a mapping from $H \backslash G$ to itself by
sending a right coset $H \, y$ to $H \, y \, g^{-1}$.  In this way we
get homomorphisms from $G$ into the groups of permutations on $G / H$
and $H \backslash G$.

	Notice that each left coset $x \, H$ and right coset $H \, y$
have the same number of elements as $H$ does, for all $x, y \in G$.
It follows that the number of elements in $G$ is equal to the number
of elements in $H$ times the number of elements in $G / H$ or $H
\backslash G$, which have the same number of elements.

	If $H$ is a normal subgroup of $G$, then every left coset of
$G$ is a right coset of $G$, and vice-versa.  Standard arguments show
that the group operation on $G$ induces a group operation on the
quotient space $G / H$ in a natural way, so that the quotient mapping
from $G$ onto $G / H$ is a group homomorphism.  Thus every normal
subgroup of $G$ is the kernel of a homomorphism from $G$ into some
group.

	As a basic example, the integers ${\bf Z}$ form an abelian
group under addition.  For each positive integer $m$, the set $m \,
{\bf Z}$ of integer multiples of $m$ is a subgroup of ${\bf Z}$.

	Because ${\bf Z}$ is abelian, $m \, {\bf Z}$ is automatically
a normal subgroup of ${\bf Z}$.  The quotient ${\bf Z} / m \, {\bf Z}$
is the group of integers modulo $m$ under addition, which has $m$
elements

	If $G$ is any group, then $Z(G)$ denotes the center of $G$,
consisting of those elements of $G$ which commute with all other
elements of $G$, i.e.,
\begin{equation}
	Z(G) = \{x \in G : x \, y = y \, x \hbox{ for all } y \in G\}.
\end{equation}
Clearly $Z(G)$ is a normal subgroup of $G$.  Moreover, $Z(G)$ is
invariant under any automorphism of $G$.

	If $G$ is a group, and $x$, $y$ are elements of $G$, then
the commutator of $x$ and $y$ is given by
\begin{equation}
	x \, y \, x^{-1} \, y^{-1}.
\end{equation}
The set of commutators in $G$ is invariant under taking inverses and
under arbitrary automorphisms of $G$.  The commutator subgroup of $G$
is the subgroup generated by the commutators, consisting of all finite
products of commutators, and it is a normal subgroup.  Indeed, it is
invariant under automorphisms of $G$ by construction.  The quotient of
$G$ by the commutator subgroup is an abelian group, and the kernel of
any homomorphism from $G$ into an abelian group contains the
commutator subgroup.

\section{Fields, vector spaces}
\setcounter{equation}{0}

	Let $k$ be a field.  Thus $k$ is a nonempty set with two
distinguished elements $0$, $1$, $0 \ne 1$, and equipped with binary
operations of addition $+$ and multiplication $\cdot$.  These
operations satisfy the usual associative, commutative, and
distributive laws, $0$, $1$ are the additive and multiplicative
identity elements for $k$, respectively, each element $x$ of $k$ has
an additive inverse $-x$, which implies that $k$ is a commutative
group with respect to addition, and each nonzero element $x$ of $k$
has a multiplicative inverse $x^{-1}$, which implies that the nonzero
elements of $k$ form a commutative group with respect to
multiplication.

	Of course the rational numbers ${\bf Q}$, the real numbers
${\bf R}$, and the complex numbers ${\bf C}$ are fields with respect
to the usual operations of addition and multiplication.  If $p$ is a
positive integer which is a prime number, which means that $1$, $p$
are the only positive integers of which $p$ is an integer multiple,
then the integers modulo $p$ form a field with respect to addition and
multiplication of integers modulo $p$.

	In general, if for each positive integer $n$ the sum of $n$
$1$'s in a field $k$ is not equal to $0$, then we say that $k$ has
characteristic $0$.  In this event $k$ contains an isomorphic copy of
the rational numbers, in the sense that there is a one-to-one mapping
from ${\bf Q}$ into $k$ which map the additive and multiplicative
identity elements to themselves and which preserve the field
operations of addition and multiplication.

	Otherwise there is a positive integer $n$ such that the sum of
$n$ $1$'s in $k$ is equal to $0$.  The smallest such positive integer
$n$ is a prime number, and is called the characteristic of the field
$k$.  If $k$ is a field with characteristic $p$, then $k$ contains a
copy of the integers modulo $p$, consisting of $0$, $1$, and the other
elements of $k$ by adding $1$'s.

	Let $k$ be a field.  A vector space over $k$ is an abelian
group $V$ with group operation $+$ and identity element $0$ which is
also equipped with an operation of scalar multiplication which permits
one to multiply an element of $V$ by an element of $k$, with the usual
compatibility conditions between addition and scalar multiplication.

	Suppose that $V$ is a vector space over a field $k$ and that
$v_1, \ldots, v_n$ are elements of $V$.  We say that $v_1, \ldots,
v_m$ are linearly dependent if there are scalars $\alpha_1, \ldots,
\alpha_m \in k$, at least one of which is different from $0$, such
that
\begin{equation}
	\alpha_1 \, v_1 + \cdots + \alpha_m \, v_m = 0.
\end{equation}
If this does not happen, then we say that $v_1, \ldots, v_m$ are
linearly independent in $V$.

	We say that vectors $w_1, \ldots, w_n$ in $V$ span $V$ if
every element of $V$ can be expressed as a linear combination of the
$w_j$'s, i.e., as a sum of scalar multiples of the $w_j$'s.
If there is a finite collection of vectors in $V$ which spans $V$,
then $V$ is said to have finite dimension.

	A finite collection of vectors in $V$ is a basis for $V$ if
the vectors are linearly independent and span $V$.  This is equivalent
to saying that every element of $V$ can be expressed in a unique way
as a linear combination of the vectors in the basis.

	Suppose that $V$ is a finite-dimensional vector space over
$k$, and that $V$ is spanned by the vectors $w_1, \ldots, w_n$.  If
$w_1, \ldots, w_n$ are linearly independent, then we have a basis for
$V$.  Otherwise, one of the $w_j$'s can be expressed as a linear
combination of the others, and we can remove such a vector and still
have a collection of vectors which spans $V$.  By repeating the
process we get a basis for $V$.  Similarly, we can start with a
collection of linearly independent vectors in a finite-dimensional
vector space and add finitely many vectors to the collection if
necessary to get a basis.

	If $V$ is a vector space over $k$ which contains a collection
of $m$ linearly independent vectors and a collection of $n$ vectors
which spans $V$, then $m \le n$.  This can be derived from the fact
that a set of $n$ homogeneous linear equations in $m$ variables has a
nontrivial solution when $m > n$.

	The dimension of a finite-dimensional vector space $V$ can be
defined as the smallest number of vectors needed to span $V$, which is
the same as the maximal number of linearly independent vectors in $V$.
If $v_1, \ldots, v_n$ is a basis for $V$, then the dimension of $V$ is
equal to $n$.  If $V$ contains only the zero vector, then $V$ has
dimension $0$, and one can interpret the empty set of vectors as being
a basis for $V$.

	If $k$ is a field and $n$ is a positive integer, then we get a
vector space $k^n$ by considering the set of $n$-tuples $x = (x_1,
\ldots, x_n)$ with each $x_j \in k$ and using coordinatewise addition
and scalar multiplication.  The standard basis in $k^n$ consists of
the $n$ vectors $e_1, \ldots, e_n$ where $e_j$ has $j$th coordinate
equal to $1$ and the rest equal to $0$.

	Suppose that $V_1$, $V_2$ are vector spaces over the same
field $k$.  A mapping $f$ from $V_1$ to $V_2$ is linear if for each
$\alpha, \beta \in k$ and $v, w \in V_1$ we have that
\begin{equation}
	f(\alpha \, v + \beta \, w) = \alpha \, f(v) + \beta \, f(w).
\end{equation}
In particular $f$ maps the zero element of $V_1$ to the zero element
of $V_2$.

	The space of linear mappings from $V_1$ to $V_2$ is denoted
$\mathcal{L}(V_1, V_2)$.  Linear mappings from $V_1$ to $V_2$ can be
added or multiplied by elements of $k$ to get linear mappings again,
so that $\mathcal{L}(V_1, V_2)$ is a vector space over $k$ in a
natural way.

	If $V_1$, $V_2$, and $V_3$ are vector spaces over $k$ and $f_1
: V_1 \to V_2$, $f_2 : V_2 \to V_3$ are linear mappings, then the
composition $f_2 \circ f_1$, which is the mapping from $V_1$ to $V_3$
given by
\begin{equation}
	(f_2 \circ f_1)(v) = f_2(f_1(v))
\end{equation}
for $v \in V_1$, is linear as well.  If $f$ is a linear mapping from
$V_1$ to $V_2$ which is a one-to-one mapping of $V_1$ onto $V_2$, then
the inverse mapping from $V_2$ to $V_1$ is linear too, and we say that
$f$ is an isomorphism from $V_1$ onto $V_2$.

	Suppose that $V$ is a vector space over a field $k$, $n$ is a
positive integer, and $v_1, \ldots, v_n$ are elements of $V$.  There
is a unique linear mapping from $k^n$ into $V$ which takes the
standard basis vector $e_j$ to $v_j$ for each $j$, $1 \le j \le n$.
The vectors $v_1, \ldots, v_n$ in $V$ are linearly independent if and
only if the linear mapping is one-to-one.  The span of $v_1, \ldots,
v_n$ is equal to $V$ if and only if $f$ maps $k^n$ onto $V$.  In
particular, $v_1, \ldots, v_n$ form a basis of $V$ if and only if this
linear mapping is an isomorphism of $k^n$ onto $V$.

	If $V$ is a vector space over $k$, then a linear subspace of
$V$ is a subset of $V$ containing $0$ and which is closed under
addition of vectors and scalar multiplication.  In other words, a
linear subspace of $V$ is a vector space over $k$ too using the
restriction of the vector space operations from $V$.  If $V$ has
finite dimension, then every linear subspace of $V$ has finite
dimension less than or equal to the dimension of $V$.

	If $V_1$, $V_2$ are vector spaces over $k$ and $f$ is a linear
mapping from $V_1$ into $V_2$, then the image of $V_1$ in $V_2$ under
$f$ is a linear subspace of $V_2$.  The kernel of $f$ is the set of
vectors $v \in V_1$ such that $f(v) = 0$ in $V_1$, and it is a linear
subspace of $V_1$.  If $V_1$ is finite-dimensional, then the image of
$V_1$ in $V_2$ is finite-dimensional, with dimension less than or
equal to the dimension of $V_1$.  More precisely, the dimension of
$V_1$ is equal to the sum of the dimensions of the kernel of $f$ and
the image of $f$.  One can show this by combining a basis for the
kernel of $f$ with a collection of vectors in $V_1$ whose images form
a basis for the image of $V_1$ under $f$ to get a basis for $V_1$.

	If $V$ is a vector space over $k$ and $W$ is a linear subspace
of $V$, then we can form the quotient $V / W$.  One can think of this
first as a quotient of abelian groups with respect to addition, and
then check that scalar multiplication is well-defined on the quotient.
There is a canonical quotient mapping which is a linear mapping from
$V$ onto $V / W$ whose kernel is equal to $W$.

	Suppose that $V_1$, $V_2$ are finite-dimensional vector spaces
over $k$.  The vector space $\mathcal{L}(V_1, V_2)$ of linear mappings
from $V_1$ to $V_2$ is a finite-dimensional vector space over $k$ with
dimension equal to the product of the dimensions of $V_1$, $V_2$.  To
see this one can choose bases for $V_1$, $V_2$ and characterize linear
transformations from $V_1$ into $V_2$ by a matrix of coefficients
which specify how a basis vector in $V_1$ is mapped to a linear
combination of basis vectors in $V_2$.

	If $k$ is a finite field with $N$ elements and $V$ is a
finite-dimensional vector space over $k$ with dimension equal to $n$,
then $V$ has $N^n$ elements.  This follows from the fact that there is
an isomorphism between $V$ and $k^n$.

	If $k$ is a finite field with characteristic $p$, then we can
view $k$ as a vector space over the integers modulo $p$.  If $k$ has
finitely many elements, then the number $N$ of elements of $k$ is
equal to $p^l$ for some positive integer $l$, by the remarks of the
preceding paragraph.

	For that matter, any field can be viewed as a vector space
over a subfield.

\section{Algebras}
\setcounter{equation}{0}

	Let $k$ be a field.  To say that $\mathcal{A}$ is an algebra
over $k$ means that $\mathcal{A}$ is a vector space over $k$ equipped
with a binary operation
\begin{equation}
	(a, b) \mapsto a \, b
\end{equation}
which is linear in each of $a$, $b$.  More precisely, for each $a \in
\mathcal{A}$, the mapping $b \mapsto a \, b$ should be a linear
mapping from $\mathcal{A}$ into itself, and for each $b \in
\mathcal{A}$ the mapping $a \mapsto a \, b$ should be a linear mapping
from $\mathcal{A}$ into itself.

	If
\begin{equation}
	(a \, b) \, c = a \, (b \, c)
\end{equation}
for all $a, b, c \in \mathcal{A}$, then we say that $\mathcal{A}$ is
an associative algebra over $k$.  Sometimes this is included in the
definition of an algebra, but it will be convenient here to mention it
separately.

	If $e$ is an element of $\mathcal{A}$ such that
\begin{equation}
	a \, e = e \, a = a
\end{equation}
for all $a \in \mathcal{A}$, then we say that $e$ is a multiplicative
identity element in $\mathcal{A}$.  Clearly $e$ is unique when it
exists.

	If
\begin{equation}
	a \, b = b \, a
\end{equation}
for all $a, b \in \mathcal{A}$, then we say that $\mathcal{A}$
is a commutative algebra.
	
	If $E$ is a nonempty set, then the space of $k$-valued
functions on $E$ is a commutative algebra over $k$ with respect to
ordinary addition and multiplication of functions, with multiplicative
identity element given by the constant function equal to $1$ at every
point in $E$.  If $E$ is a finite set, then this algebra is
finite-dimensional as a vector space over $k$, with dimension as a
vector space equal to the number of elements of $E$.

	Let $V$ be a vector space over $k$, and let $\mathcal{L}(V)$
denote the space of linear transformations from $V$ to itself.  This
is an algebra over $k$ using ordinary addition and scalar
multiplication of linear transformations and composition of linear
operators on $V$ as multiplication, and the identity transformation
$I$ on $V$, which takes every element of $V$ to itself, is the
multiplicative identity element in $\mathcal{L}(V)$.  If $V$ is a
finite-dimensional vector space over $k$, then $\mathcal{L}(V)$ is
finite-dimensional as a vector space over $k$, with dimension equal to
the square of the dimension of $V$.

	Suppose that $\mathcal{A}$ is an algebra over $k$, and that
$\mathcal{A}$ is finite-dimensional as a vector space over $k$.  Let
$a_1, \ldots, a_n$ be a basis for $\mathcal{A}$ as a vector space over
$k$.  We can write the product $a_j \, a_l$ as a linear combination of
the basis vectors for each $j$, $l$.  This leads to a family of $n^3$
elements of $k$ which describe multiplication in $\mathcal{A}$, and
any choice of $n^3$ coefficients in $k$ leads to an algebra structure
on an $n$-dimensional vector space over $k$ with a specified basis.

	In some situations one considers vector spaces or algebras
with additional structure.  This may involve distinguished subspaces,
as with filtrations or gradings.  There might be some extra operators
present.  There could be something like a topological structure which
is helpful, perhaps related to some kind of norm or family of norms.
At any rate, it is often easy and interesting to accommodate
additional ingredients like these into the basic notions.

\section{Lie algebras}
\setcounter{equation}{0}

	Let $k$ be a field.  To say that $\lambda$ is a Lie algebra
over $k$ means that $\lambda$ is a vector space over $k$ equipped with
a binary operation
\begin{equation}
	(x, y) \mapsto [x, y]
\end{equation}
which is linear in $x, y \in \lambda$ and satisfies 
\begin{equation}
	[x, x] = 0
\end{equation}
for all $x \in \lambda$ and the Jacobi identity
\begin{equation}
	[x, [y, z]] + [y, [z, x]] + [z, [x, y]] = 0
\end{equation}
for all $x, y, z \in \mathcal{\lambda}$.  If $x, y \in \lambda$, then
\begin{equation}
	[x + y, x + y] = [x, x] + [x, y] + [y, x] + [y, y]
\end{equation}
and hence
\begin{equation}
	[x, y] = - [y, x].
\end{equation}

	Thus a Lie algebra is an algebra in the general sense
described in the previous section.  A Lie algebra is said to be
commutative if the bracket of any two elements of the algebra is equal
to $0$.  The general definition of commutativity of an algebra over
$k$ would say that the bracket should be symmetric, which is
equivalent to saying that it is identically equal to $0$ if $k$ does
not have characteristic equal to $2$.  If $k$ has characteristic equal
to $2$, then the bracket is automatically symmetric, and commutativity
as a Lie algebra means that it is equal to $0$.  Notice that if $V$ is
any vector space over $k$, then one can define a Lie bracket on $V$ by
saying that the bracket of any two vectors is equal to $0$.

	If $\mathcal{A}$ is an associative algebra over $k$, then we
can define
\begin{equation}
	[a, b] = a \, b - b \, a
\end{equation}
for all $a, b \in \mathcal{A}$.  One can check that this satisfies the
conditions of a Lie algebra, using associativity of the product in
$\mathcal{A}$ to get the Jacobi identity for this bracket.  We shall
sometimes write $\lambda(\mathcal{A})$ for $\mathcal{A}$ as a Lie
algebra with this bracket.  In particular, we can define a Lie bracket
on the linear transformations on a vector space $V$ over $k$ by
\begin{equation}
	[A, B] = A \, B - B \, A = A \circ B - B \circ A
\end{equation}
for all $A, B \in \mathcal{L}(V)$.  Let us write $\Lambda(V)$ for the
space of linear transformations on $V$, as a Lie algebra using this
bracket, which is the same as $\lambda(\mathcal{L}(V))$.

\section{Subalgebras, ideals, homomorphisms}
\setcounter{equation}{0}

	Many familiar notions can be formulated for general algebras
over a field $k$ and their basic properties verified in the usual way.
For instance, a subalgebra of an algebra is a vector subspace which is
closed under multiplication, and a homomorphism between two algebras
is a linear mapping which preserves multiplication.

	Let $\mathcal{A}$ be an algebra over $k$, and let
$\mathcal{I}$ be a subalgebra of $\mathcal{A}$.  We say that
$\mathcal{I}$ is a left ideal in $\mathcal{A}$ if for each $x \in
\mathcal{A}$ and $y \in \mathcal{I}$ we have that $x \, y \in
\mathcal{I}$.  We say that $\mathcal{I}$ is a right ideal in
$\mathcal{A}$ if for each $y \in \mathcal{I}$ and $z \in \mathcal{A}$
we have that $y \, z \in \mathcal{I}$.  If $\mathcal{I}$ is both a
left and right ideal in $\mathcal{A}$, then we say that $\mathcal{I}$
is a two-sided ideal in $\mathcal{A}$.  For a commutative algebra
$\mathcal{A}$ these three notions coincide.

	If $\lambda$ is a Lie algebra over $k$, then a subalgebra of
$\lambda$ is a linear subspace of $\lambda$ which is closed under
brackets.  A subalgebra $\iota$ of $\lambda$ is an ideal if for each
$x \in \lambda$ and $y \in \iota$ we have that $[x, y] \in \iota$.  As
for a commutative algebra, we do not need to discriminate between
left, right, and two-sided ideals in a Lie algebra.

	Let $\mathcal{A}_1$, $\mathcal{A}_2$ be algebras over $k$, and
suppose that $h$ is a homomorphism from $\mathcal{A}_1$ into
$\mathcal{A}_2$.  By definition, the kernel of $h$ is the set of $x
\in \mathcal{A}_1$ such that $h(x) = 0$ in $\mathcal{A}_2$, which is
automatically a two-sided ideal in $\mathcal{A}_1$.  Conversely,
suppose that $\mathcal{A}$ is an algebra over $k$ and that
$\mathcal{I}$ is a two-sided ideal in $\mathcal{A}$.  One can define
the quotient $\mathcal{A} / \mathcal{I}$ as an algebra over $k$ with a
canonical quotient homomorphism from $\mathcal{A}$ onto $\mathcal{A} /
\mathcal{I}$ with kernel equal to $\mathcal{I}$.  Specifically, one
can define $\mathcal{A} / \mathcal{I}$ first as a vector space over
$k$, and then check that multiplication is well-defined on the
quotient.

	If $\mathcal{A}$ is an associative algebra over $k$, then
define the center $Z(\mathcal{A})$ to be the set of $x \in
\mathcal{A}$ such that $x \, y = y \, x$ for all $y \in \mathcal{A}$.
One can check that this is a subalgebra of $\mathcal{A}$ which is
invariant under automorphisms of $\mathcal{A}$.  If $\mathcal{A}$ has
a multiplicative identity element $e$, then $e \in Z(\mathcal{A})$.
In general $Z(\mathcal{A})$ is not an ideal in $\mathcal{A}$, but it
is automatically a two-sided ideal when it is a one-sided ideal.

	If $\lambda$ is a Lie algebra over $k$, then the center
$Z(\lambda)$ of $\lambda$ is the set of $x \in \lambda$ such that $[x,
y] = 0$ for all $y \in \lambda$.  This is an ideal in $\lambda$.  The
center of $\lambda$ is invariant under automorphisms of $\lambda$ too.

	By a representation of an associative algebra $\mathcal{A}$
over $k$ on a vector space $V$ over $k$ we mean a homomorphism of
$\mathcal{A}$ into the algebra $\mathcal{L}(V)$ of linear
transformations on $V$.  By a representation of a Lie algebra
$\lambda$ over $k$ on a vector space $V$ over $k$ we mean a
homomorphism from $\lambda$ into the Lie algebra $\Lambda(V)$ of
linear transformations on $V$.

	Let $\mathcal{A}$ be an associative algebra over $k$.  We can
define a representation of $\mathcal{A}$ on itself, as a vector over
$k$, by sending $a \in \mathcal{A}$ to the linear transformation
\begin{equation}
	x \mapsto a \, x
\end{equation}
on $\mathcal{A}$.  The kernel of this homomorphism is the two-sided
ideal consisting of $a \in \mathcal{A}$ such that $a \, x = 0$ for all
$x \in \mathcal{A}$.

	Now let $\lambda$ be a Lie algebra over $k$.  For each $x \in
\lambda$, define $\ad_x$ as a mapping from $\lambda$ to itself by
\begin{equation}
	{\ad}_x(y) = [x, y].
\end{equation}
A straightforward computation using the Jacobi identity shows that
\begin{equation}
	x \mapsto {\ad}_x
\end{equation}
defines a representation of $\lambda$ on $\lambda$ as a vector space
over $k$, the adjoint representation of $\lambda$.  The kernel of this
representation consists exactly of the center of $\lambda$.

	If $\mathcal{A}$ is an associative algebra over $k$, then the
commutator subalgebra of $\mathcal{A}$ is the subalgebra of
$\mathcal{A}$ generated by commutators, i.e., consisting of finite
linear combinations of products of commutators $x \, y - y \, x$, $x,
y \in \mathcal{A}$.  This is a subalgebra of $\mathcal{A}$ by
construction which is invariant under automorphisms, which is not an
ideal in general, and which is a two-sided ideal when it is a
one-sided ideal.  If $\lambda$ is a Lie algebra over $k$, then the set
of linear combinations of brackets $[x, y]$, $x, y \in \lambda$, is an
ideal in $\lambda$ which is invariant under automorphisms of
$\lambda$.  The quotient of $\lambda$ by its commutator ideal is
automatically a commutative Lie algebra.  Of course in any algebra
$\mathcal{A}$ over $k$ one can consider the two-sided ideal consisting
of linear combinations of arbitrary products, and the quotient of
$\mathcal{A}$ by this ideal is an algebra in which every product is
equal to $0$.

\section{Derivations}
\setcounter{equation}{0}

	If $\mathcal{A}$ is an algebra over a field $k$, then a
derivation $\delta$ on $\mathcal{A}$ is a linear mapping from 
$\mathcal{A}$ to itself such that
\begin{equation}
	\delta(a \, b) = \delta(a) \, b + a \, \delta(b)
\end{equation}
for all $a, b \in \mathcal{A}$.  Linear combinations of derivations on
$\mathcal{A}$ with coefficients in $k$ are derivations on
$\mathcal{A}$, and if $\delta_1$, $\delta_2$ are derivations on
$\mathcal{A}$, then
\begin{equation}
	[\delta_1, \delta_2] 
		= \delta_1 \circ \delta_2 - \delta_2 \circ \delta_1
\end{equation}
is a derivation on $\mathcal{A}$.  Thus the derivations on
$\mathcal{A}$ form a Lie algebra, a subalgebra of the Lie algebra of
all linear transformations on $\mathcal{A}$ with respect to the usual
bracket of linear transformations.  In particular this applies to
associative algebras and to Lie algebras.  One can check that a
derivation on an associative algebra $\mathcal{A}$ is also a
derivation on the corresponding Lie algebra $\lambda(\mathcal{A})$.

	If $\lambda$ is a Lie algebra over $k$, then a derivation on
$\lambda$ is a linear mapping $\delta$ of $\lambda$ into itself such
that $\delta([x, y]) = [\delta(x), y] + [x, \delta(y)]$ for all $x, y
\in \lambda$.  A straightforward computation using the Jacobi identity
shows that for each $x \in \lambda$, $\ad_x$ is a derivation on
$\lambda$, which means that the adjoint representation of a Lie
algebra $\lambda$ maps into the Lie algebra of derivations on
$\lambda$.  If $\mathcal{A}$ is an associative algebra over $k$, then
one can check that the adjoint representation of
$\lambda(\mathcal{A})$ maps into the Lie algebra of derivations of
$\mathcal{A}$.

	Suppose now that $\mathcal{A}$ is a commutative associative
algebra over $k$.  If $a$ is an element of $\mathcal{A}$ and $\delta$
is a derivation on $\mathcal{A}$, then $a \, \delta$ defines a
derivation on $\mathcal{A}$, where $(a \, \delta)(x) = a \, \delta(x)$
for $x \in \mathcal{A}$.

	Let $\delta_1, \ldots, \delta_n$ be derivations on
$\mathcal{A}$ which commute as operators on $\mathcal{A}$, and let
$\mathcal{A}^n$ denote the space of $n$-tuples of elements of
$\mathcal{A}$, which is a vector space over $k$, using coordinatewise
addition and scalar multiplication.  For each $a = (a_1, \ldots,
a_n)$, $b = (b_1, \ldots, b_n)$ in $\mathcal{A}^n$, let $[a, b]$ be
the element of $\mathcal{A}^n$ whose $j$th component is equal to
\begin{equation}
	\sum_{l=1}^n a_l \, \delta_l(b_j) - b_l \, \delta_l(a_j).
\end{equation}
One can check that $\mathcal{A}^n$ becomes a Lie algebra with this
choice of bracket.  For each $a = (a_1, \ldots, a_n) \in
\mathcal{A}^n$,
\begin{equation}
	x \in \mathcal{A} \mapsto \sum_{j=1}^n a_j \, \delta_j(x)
\end{equation}
defines a derivation on $\mathcal{A}$.  The bracket on $\mathcal{A}^n$
just described corresponds exactly to the commutator of the associated
derivations on $\mathcal{A}$, by construction.

\section{Smooth functions on ${\bf R}^n$}
\setcounter{equation}{0}

	Fix a positive integer $n$, and let $C^\infty({\bf R}^n)$
denote the space of smooth real-valued functions on ${\bf R}^n$, i.e.,
the space of continuous real-valued functions $f(x)$ on ${\bf R}^n$
such that the partial derivatives of $f(x)$ of all orders exist and
are continuous.  As usual, the sum and product of smooth functions is
again a smooth function.  Thus $C^\infty({\bf R}^n)$ is a commutative
algebra over the real numbers.  It is customary to equip
$C^\infty({\bf R}^n)$ with a topology such that a sequence of smooth
functions $\{f_j\}_{j=1}^\infty$ on ${\bf R}^n$ converges to a smooth
function $f$ on ${\bf R}^n$ if the $f_j$'s converge to $f$ uniformly
on compact subsets of ${\bf R}^n$, and if the derivatives of the
$f_j$'s converge to the corresponding derivatives of $f$ uniformly on
compact subsets of ${\bf R}^n$.  Standard arguments show that sums and
products of convergent sequences of smooth functions converge to the
corresponding sums and products of the limits of the sequences.

	By a smooth vector field on ${\bf R}^n$ we mean an ${\bf
R}^n$-valued function $V(x) = (V_1(x), \ldots, V_n(x))$ on ${\bf R}^n$
whose components $V_1(x), \ldots, V_n(x)$ are smooth functions on
${\bf R}^n$.  If $f(x)$ is a real-valued smooth function on ${\bf
R}^n$ and $V(x)$ is a smooth vector field on ${\bf R}^n$, then let
$V(f)$ be the smooth function on ${\bf R}^n$ defined by
\begin{equation}
  (V(f))(x) = \sum_{j=1}^n V_j(x) \, \frac{\partial}{\partial x_j} \, f(x),
\end{equation}
which is the same as the directional derivative of $f$ in the
direction of $V$ at $x$.  Thus $f \mapsto V(f)$ is a continuous linear
mapping on $C^\infty({\bf R}^n)$ which is a derivation, since
\begin{equation}
	V(f_1 \, f_2) = V(f_1) \, f_2  + f_1 \, V(f_2)
\end{equation}
for all $f_1, f_2 \in C^\infty({\bf R}^n)$ by the usual Leibniz rule
from calculus.  

	If $V = (V_1(x), \ldots, V_n(x))$, $W = (W_1(x), \ldots, W_n)$
are smooth vector fields on ${\bf R}^n$, then their Lie bracket $[V,
W]$ is defined to be the smooth vector field on ${\bf R}^n$ whose
$j$th component is equal to
\begin{equation}
	\sum_{l=1}^n V_l(x) \, \frac{\partial}{\partial x_l} \, W_j(x)
		- W_l(x) \, \frac{\partial}{\partial x_l} \, V_j(x).
\end{equation}
For each smooth function $f$ on ${\bf R}^n$ we have that
\begin{equation}
	[V, W](f) = V(W(f)) - W(V(f)),
\end{equation}
which says that the derivation on $C^\infty({\bf R}^n)$ associated to
$[V, W]$ is the commutator of the derivations associated to $V$ and
$W$.

\section{Polynomials}
\setcounter{equation}{0}

	Let $k$ be a field, $\mathcal{A}$ an associative algebra over
$k$, and $n$ be a positive integer.  Let $\mathcal{A}[t_1, \ldots,
t_n]$ denote the usual polynomial algebra with coefficients in
$\mathcal{A}$ in the indeterminants $t_1, \ldots, t_n$.

	More precisely, recall that a multi-index is an $n$-tuple
$\alpha = (\alpha_1, \ldots, \alpha_n)$ of nonnegative integers.  We
can add multi-indices coordinatewise, and the degree of a multi-index
$\alpha$ is defined to be the sum of its coordinates.  If $\alpha$ is
a multi-index, then $t^\alpha$ denotes the monomial
\begin{equation}
	t^\alpha = t_1^{\alpha_1} \cdots t_n^{\alpha_n},
\end{equation}
where $t_j^{\alpha_j}$ is interpreted as being equal to $1$ when
$\alpha_j = 0$, and $t^\alpha$ is interpreted as being equal to $1$
when $\alpha = 0$.  An element of $\mathcal{A}[t_1, \ldots, t_n]$ can
be expressed as
\begin{equation}
	\sum_{\alpha_1 + \cdots \alpha_n \le m} a_\alpha \, t^\alpha,
\end{equation}
where $m$ is a nonnegative integer, the $c_\alpha$'s are elements of
$\mathcal{A}$, and the sum is taken over multi-indices $\alpha$ with
degree less than or equal to $m$.  Thus $\mathcal{A}[t_1, \ldots,
t_n]$ contains a copy of $\mathcal{A}$ as the constant polynomials,
i.e., the polynomials with only a constant term.

	We can add polynomials and multiply them by elements of $k$
termwise.  We can also multiply polynomials, where
\begin{equation}
	t^\alpha \, t^\beta = t^{\alpha + \beta}
\end{equation}
for arbitrary multi-indices $\alpha$, $\beta$, and where the monomials
$t^\alpha$ commute with each other and with elements of $\mathcal{A}$
by definition. In this way $\mathcal{A}[t_1, \ldots, t_n]$ becomes an
associative algebra over $k$.

	For each nonnegative integer $\ell$, let us write
$\mathcal{A}_\ell[t_1, \ldots, t_n]$ for the polynomials with
coefficients in $\mathcal{A}$ which are homogeneous of degree equal to
$\ell$, which is to say polynomials of the form
\begin{equation}
	\sum_{\alpha_1 + \cdots + \alpha_n = \ell} a_\alpha \, t^\alpha.
\end{equation}
Thus $\mathcal{A}_\ell[t_1, \ldots, t_n]$ is a vector subspace of
$\mathcal{A}[t_1, \ldots, t_n]$ as a vector space over $k$ which is
closed under multiplication by elements of $\mathcal{A}$.

	For each $j$, $1 \le j \le n$, we can define the operator
$\partial_j$ on $\mathcal{A}[t_1, \ldots, t_n]$ in the usual way, by
formally differentiating in $t_j$.  Specifically, if $\alpha$ is a
multi-index and $1 \le j \le n$, let $r_j(\alpha)$ be the multi-index
which agrees with $\alpha$ except for the $j$th component, where the
$j$th component of $r_j(\alpha)$ is equal to $\alpha_j - 1$ when
$\alpha_j \ge 1$ and to $0$ when $\alpha_j = 0$.  We define
$\partial_j$ acting on polynomials with coefficients in $\mathcal{A}$
by
\begin{equation}
	\partial_j \, \bigg(\sum_{|\alpha| \le m} a_\alpha \, t^\alpha \bigg)
	= \sum_{|\alpha| \le m} \alpha_j \, a_\alpha \, t^{r_j(\alpha)}.
\end{equation}
Clearly $\partial_j$ is linear on $\mathcal{A}[t_1, \ldots, t_n]$ as a
vector space over $k$, and it is linear with respect to multiplication
on the left or right by elements of $\mathcal{A}$.  For each $j$,
$\partial_j$ is a derivation on the algebra $\mathcal{A}[t_1, \ldots,
t_n]$, and $\partial_j$ commutes with $\partial_l$ for all $j$, $l$.

	Now suppose that $\mathcal{A}$ is a commutative associative
algebra over $k$.  For each positive integer $n$, the polynomial
algebra $\mathcal{A}[t_1, \ldots, t_n]$ is a commutative associative
algebra over $k$ too.

	Let $p_1, \ldots, p_n$ be elements of $\mathcal{A}[t_1,
\ldots, t_n]$, and consider the operator on $\mathcal{A}[t_1, \ldots,
t_n]$ given by
\begin{equation}
	f \in \mathcal{A}[t_1, \ldots, t_n]
		\mapsto \sum_{j=1}^n p_j \, \partial_j (f).
\end{equation}
This operator is linear on $\mathcal{A}[t_1, \ldots, t_n]$ as a vector
space over $k$, and also linear with respect to multiplication by
elements of $\mathcal{A}$.  Moreover it defines a derivation on
$\mathcal{A}[t_1, \ldots, t_n]$.

	Suppose that $q_1, \ldots, q_n$ are also elements of
$\mathcal{A}[t_1, \ldots, t_n]$.  Define $r_1, \ldots, r_n$ in
$\mathcal{A}[t_1, \ldots, t_n]$ by
\begin{equation}
	r_j = \sum_{l = 1}^n p_l \, \partial_l (q_j) - q_l \, \partial_l(p_j).
\end{equation}
This defines a bracket on the space of $n$-tuples of polynomials with
coefficients in $\mathcal{A}$ so that it becomes a Lie algebra over
$k$.  By construction, the differential operator $\sum_{j=1}^n r_j \,
\partial_j$ is equal to the commutator of the operators $\sum_{l=1}^n
p_l \, \partial_l$ and $\sum_{m=1}^n q_m \, \partial_m$ on
$\mathcal{A}[t_1, \ldots, t_n]$.  In other words, we get a Lie bracket
on $n$-tuples of polynomials which corresponds exactly to the
commutator of the associated first-order differential operators on
$\mathcal{A}[t_1, \ldots, t_n]$.

	If $p_1, \ldots, p_n$ are elements of $\mathcal{A}_1[t_1,
\ldots, t_n]$, then the differential operator $\sum_{j=1}^n p_j \,
\partial_j$ maps $\mathcal{A}_\ell[t_1, \ldots, t_n]$ to itself for
each nonnegative integer $\ell$.  Moreover, if $q_1, \ldots, q_n$ are
elements of $\mathcal{A}_1[t_1, \ldots, t_n]$, then $r_1, \ldots, r_n$
as defined in the preceding paragraph are contained in
$\mathcal{A}_1[t_1, \ldots, t_n]$.  Thus $n$-tuples of homogeneous
polynomials of degree $1$ with coefficients in $\mathcal{A}$ form a
Lie subalgebra of the Lie algebra $n$-tuples of polynomials with
coefficients in $\mathcal{A}$ described in the previous paragraph.

\section{Matrices}
\setcounter{equation}{0}

	Let $k$ be a field and let $\mathcal{A}$ be an associative
algebra over $k$.  Fix a positive integer $n$.

	Let us write $\mathcal{A}^n$ for the space of $n$-tuples $x =
(x_1, \ldots, x_n)$ with $x_j \in \mathcal{A}$ for $1 \le j \le n$.
This is a vector space over $k$ with respect to coordinatewise
addition and scalar multiplication by $k$.  If $\mathcal{A}$ is
finite-dimensional as a vector space over $k$, then $\mathcal{A}^n$ is
a finite-dimensional vector space over $k$, with dimension equal to
$n$ times the dimension of $\mathcal{A}$.

	Let us write $M_n(\mathcal{A})$ for the space of $n \times n$
matrices with entries in $\mathcal{A}$.  We can add matrices and
define scalar multiplication by $k$ entry by entry.  We can also
multiply matrices in the usual manner.  Specifically, if $(a_{j, l})$
and $(b_{p, q})$ are $n \times n$ matrices with entries in
$\mathcal{A}$, then their product is defined to be the $n \times n$
matrix
\begin{equation}
	c_{j, m} = \sum_{l=1}^n a_{j, l} \, b_{l, m}.
\end{equation}
In this way $M_n(\mathcal{A})$ becomes an associative algebra over
$k$.  Of course $M_n(\mathcal{A})$ is the same as $\mathcal{A}$ when
$n = 1$.

	Suppose that $(t_{j, l})$ is an $n \times n$ matrix with
entries in $\mathcal{A}$.  Define a mapping $T : \mathcal{A}^n \to
\mathcal{A}^n$ by $y = T(x)$,
\begin{equation}
	y_j = \sum_{l=1}^n t_{j, l} \, x_l,
\end{equation}
for each $x \in \mathcal{A}^n$.  Clearly $T$ is a linear mapping on
$\mathcal{A}^n$ as a vector space over $k$.  It is also linear with
respect to multiplication of vectors in $\mathcal{A}^n$ by elements of
$\mathcal{A}$ on the right side.

	The sum of two matrices with entries in $\mathcal{A}$
corresponds to the sum of the operators on $\mathcal{A}^n$ associated
to the two individual matrices.  The product of a matrix with an
element of $k$ corresponds to the operator on $\mathcal{A}^n$ which is
the product of the operator associated to the initial matrix and the
same scalar.  The product of two matrices corresponds to the operator
on $\mathcal{A}^n$ which is the composition of the operators
associated to the two matrices.

	Since $M_n(\mathcal{A})$ is an associative algebra over $k$,
we obtain a Lie algebra $\lambda(M_n(\mathcal{A}))$, where the bracket
of two matrices with entries in $\mathcal{A}$ is given by the
commutator defined using matrix multiplication.  This Lie algebra is
often denoted $gl(n, \mathcal{A})$.

	Let $T = (t_{j, l})$ be an $n \times n$ matrix with entries in
$\mathcal{A}$.  The trace of $T$ is denoted $\tr T$ and
defined by
\begin{equation}
	\tr T = \sum_{j=1}^n t_{j, j}.
\end{equation}
The trace is linear as a mapping from $M_n(\mathcal{A})$ into
$\mathcal{A}$ as vector spaces over $k$, and with respect to
multiplication of a matrix with entries in $\mathcal{A}$ by an element
of $\mathcal{A}$ on the left or on the right.

	Suppose that $\mathcal{A}$ is a commutative associative
algebra over $k$.  A standard computation shows that
\begin{equation}
	\tr T_1 \, T_2 = \tr T_2 \, T_1
\end{equation}
for all $T_1, T_2 \in M_n(\mathcal{A})$.  Hence the trace of any
commutator of matrices with entries in $\mathcal{A}$ is equal to $0$.
The matrices with entries in $\mathcal{A}$ and trace equal to $0$
therefore form a Lie subalgebra of $gl(n, \mathcal{A})$, which
is denoted $sl(n, \mathcal{A})$.

	Let us continue to assume to $\mathcal{A}$ is a commutative
associative algebra over $k$, and let $(a_{j, l})$ be an $n \times n$
matrix with entries in $\mathcal{A}$.  This leads to the polynomials
$p_j = \sum_{l=1}^n a_{j, l} \, t_l$ for $1 \le j \le n$, which are
homogeneous polynomials of degree $1$ with coefficients in
$\mathcal{A}$.  The commutator of two matrices with respect to matrix
multiplication corresponds to $-1$ times the bracket of $n$-tuples of
homogeneous polynomials of degree $1$ described in the previous
section.

\section{Invertible matrices}
\setcounter{equation}{0}

	Let $k$ be a field, and let $\mathcal{A}$ be an associative
algebra over $k$ with a nonzero multiplicative identity element $e$.

	In this case the algebra $M_n(\mathcal{A})$ of $n \times n$
matrices with entries in $\mathcal{A}$ also has a nonzero identity
element, namely the matrix with entries equal to $e$ on the diagonal
and to $0$ off of the diagonal.

	In the previous section we saw that each matrix in
$M_n(\mathcal{A})$ determines a mapping from $\mathcal{A}^n$ to itself
which is linear on $\mathcal{A}^n$ as a vector space over $k$, and
also linear with respect to multiplication on the right by elements of
$\mathcal{A}$.  Of course the transformation on $\mathcal{A}^n$
associated to the identity matrix in $M_n(\mathcal{A}^n)$ is the
identity transformation, which takes each element of $\mathcal{A}^n$
to itself.

	Because $\mathcal{A}$ has a multiplicative identity element,
each transformation on $\mathcal{A}^n$ which is linear on
$\mathcal{A}^n$ as a vector space over $k$ and linear with respect to
multiplication on the right by elements of $\mathcal{A}$ corresponds
to a matrix with entries in $\mathcal{A}$ in this manner.  More
precisely, any transformation of this type is determined by what it
does on the elements of $\mathcal{A}^n$ which are equal to $e$ in one
coordinate and to $0$ in the others.  The relevant matrix entries can
be obtained from the coordinates of the images of these elements of
$\mathcal{A}^n$ under the transformation.

	An element of an associative algebra with a nonzero
multiplicative identity element is said to be invertible if there is
an element of the algebra so that the product of the two in each order
is equal to the identity element in the algebra.  The invertible
elements of such an algebra form a group under multiplication.

	We can apply this to $\mathcal{A}$ or to $M_n(\mathcal{A})$.
The group of invertible elements of $M_n(\mathcal{A})$ is called the
general linear group of invertible $n \times n$ matrices with entries
in $\mathcal{A}$ and is denoted $GL(n, \mathcal{A})$.  Notice that an
element of $M_n(\mathcal{A})$ is invertible as a matrix if and only if
the corresponding transformation on $\mathcal{A}^n$ is invertible.

	Suppose that $\mathcal{A}$ is also commutative.  One can then
define a determinant function from $M_n(\mathcal{A})$ into
$\mathcal{A}$ with the usual properties.

	The determinant of the identity matrix in $M_n(\mathcal{A}^n)$
is equal to the multiplicative identity element of $\mathcal{A}$.  The
determinant of a product of matrices is equal to the corresponding
product of determinants.  As a result, the determinant of an
invertible matrix is an invertible element of $\mathcal{A}$.

	Conversely, a matrix in $M_n(\mathcal{A})$ is invertible if
its determinant is an invertible element of $\mathcal{A}$.  This
follows from standard computations in algebra, in which the product of
a matrix and another matrix determined from it is equal to the
determinant of the initial matrix times the identity matrix.

	Thus $GL(n, \mathcal{A})$ can be described as the group of $n
\times n$ matrices with entries in $\mathcal{A}$ whose determinant is
an invertible element of $\mathcal{A}$.  The determinant defines a
homomorphism from this group into the group of invertible elements of
$\mathcal{A}$.  The special linear group of $n \times n$ matrices with
entries in $\mathcal{A}$ is denoted $SL(n, \mathcal{A})$ and is the
subgroup of $M_n(\mathcal{A})$ consisting of matrices whose
determinant is equal to the identity element of $\mathcal{A}$.

\section{Real numbers}
\setcounter{equation}{0}

	If $E$ is a set of real numbers and $a$ is a real number, then
we say that $a$ is a lower bound for $E$ if $a \le x$ for all $x \in
E$.  Similarly we say that a real number $b$ is an upper bound for $E$
if $x \le b$ for all $x \in E$.

	A real number is said to be the greatest lower bound or
infimum of $E$, denoted $\inf E$, if it is a lower bound for $E$ and
if it is greater than or equal to every other lower bound of $E$.  A
real number is said to be the least upper bound or supremum of $E$,
denoted $\sup E$, if it is an upper bound for $E$ and if it is less
than or equal to every other upper bound of $E$.  It is clear from the
definitions that the infimum and supremum are unique when they exist.

	The completeness axiom for the real numbers states that every
nonempty set of real numbers which has an upper bound has a least
upper bound.  A standard consequence of this is that every nonempty
set of real numbers with a lower bound has a greatest lower bound.

	The absolute value of a real number $x$ is denoted $|x|$ and
defined to be equal to $x$ when $x \ge 0$ and to $-x$ when $x \le 0$.
The triangle inequality states that
\begin{equation}
	|x + y| \le |x| + |y|
\end{equation}
for all $x, y \in {\bf R}$.  Moreover,
\begin{equation}
	|x \, y| = |x| \, |y|
\end{equation}
for all $x, y \in {\bf R}$.

	If $\{x_j\}_{j=1}^\infty$ is a sequence of real numbers and
$x$ is a real number, then we say that $\{x_j\}_{j=1}^\infty$
converges to $x$ and write
\begin{equation}
	\lim_{j \to \infty} x_j = x
\end{equation}
if for every $\epsilon > 0$ there is an $L \ge 1$ such that
\begin{equation}
	|x - x_j| < \epsilon
\end{equation}
for all $j \ge L$.  It is easy to see that the limit is unique when it
exists.

	Suppose that $\{x_j\}_{j=1}^\infty$, $\{y_j\}_{j=1}^\infty$
are sequences of real numbers which converge to the real numbers $x$,
$y$.  It is well-known that the sequences of sums $\{x_j +
y_j\}_{j=1}^\infty$ and products $\{x_j \, y_j\}_{j=1}^\infty$
converge to the sum $x + y$ and product $x \, y$ of the limits.  If
$x_j \ne 0$ for all $j$ and $x \ne 0$, then
$\{x_j^{-1}\}_{j=1}^\infty$ converges to $x^{-1}$.

	If $\{x_j\}_{j=1}^\infty$ is a sequence of real numbers which
is monotone increasing and bounded from above, then
$\{x_j\}_{j=1}^\infty$ converges to the supremum of the $x_j$'s.
Similarly, if $\{y_j\}_{j=1}^\infty$ is a monotone decreasing sequence
of real numbers which is bounded from below, then
$\{y_j\}_{j=1}^\infty$ converges to the infimum of the $y_j$'s.

	A sequence $\{x_j\}_{j=1}^\infty$ of real numbers is a Cauchy
sequence if for each $\epsilon > 0$ there is an $L \ge 1$ such that
\begin{equation}
	|x_j - x_l| < \epsilon
\end{equation}
for all $j, l \ge L$.  It is easy to see that every convergent
sequence is a Cauchy sequence.  Conversely, another version of
completeness of the real numbers states that every Cauchy sequence of
real numbers converges.

	Indeed, suppose that $\{x_j\}_{j=1}^\infty$ is a bounded
sequence of real numbers.  For each positive integer $l$, let $w_l$ be
the infimum of $w_j$, $j \ge l$, and let $y_l$ be the supremum of
$x_j$ for $j \ge l$.  Clearly $\{w_l\}_{l=1}^\infty$ is a monotone
increasing sequence of real numbers, $\{y_l\}_{l=1}^\infty$ is a
monotone decreasing sequence of real numbers, and both sequences are
bounded.  Therefore both sequences converge, and their limits are
denoted $\liminf_{j \to \infty} x_j$ and $\limsup_{j \to \infty} x_j$,
respectively.

	Clearly $\liminf_{j \to \infty} x_j \le \limsup_{j \to \infty}
x_j$ by construction.  One can check that $\{x_j\}_{j=1}^\infty$
converges to a real number $x$ if and only if $\liminf_{j \to \infty}
x_j$ and $\limsup_{j \to \infty} x_j$ are both equal to $x$.  One can
also check that the upper and lower limits of a Cauchy sequence are
equal, which implies that the Cauchy sequence converges to their
common value.

	An infinite series $\sum_{j=1}^\infty a_j$ of real numbers is
said to converge if the sequence of partial sums $\sum_{j=1}^n a_j$
converges.  Sometimes it is convenient to begin sequences or series at
$0$ or some other integer, which works just as well for these
definitions.

	Suppose that $\sum_{j=1}^\infty a_j$ is an infinite series of
nonnegative real numbers.  The partial sums for this series are then
monotone increasing.  Such a series converges if and only if the
partial sums are bounded.

	An infinite series $\sum_{j=1}^\infty a_j$ of real numbers is
said to converge absolutely if $\sum_{j=1}^\infty |a_j|$ converges.
If a series of real numbers converges absolutely, then it converges.

	Let $p(x) = a_n \, x^n + \cdots + a_1 \, x + a_0$ be a
polynomial function on the real line with real coefficients.  If
$p(x)$ takes both positive and negative values, then one can show that
$p(x) = 0$ for some $x \in {\bf R}$.  In particular this happens if
$n$ is an odd positive integer and $a_n \ne 0$.

\section{Complex numbers}
\setcounter{equation}{0}

	A complex number $z$ can be expressed in a unique manner as $x
+ y \, i$, where $x$, $y$ are real numbers and $i$ is a specific
complex number such that $i^2 = -1$.  By definition $x$, $y$ are the
real and imaginary parts of $z$, respectively, and are denoted $\re
z$, $\im z$.

	If $z = x + y \, i$ with $x, y \in {\bf R}$, then the complex
conjugate of $z$ is denoted $\overline{z}$ and defined by
\begin{equation}
	\overline{z} = x - y \, i.
\end{equation}
For any pair of complex numbers $z$, $w$ we have that $\overline{z +
w} = \overline{z} + \overline{w}$ and $\overline{z \, w} =
\overline{z} \, \overline{w}$.

	The modulus of a complex number $z = x + y \, i$, $x, y \in
{\bf R}$, is denoted $|z|$ and defined to be $\sqrt{x^2 + y^2}$.  This
is equivalent to saying that $|z|$ is a nonnegative real number and
$|z|^2 = z \, \overline{z}$.  If $z$ happens to be a real number, then
the modulus of $z$ is the same as the absolute value of $z$.

	If $z$ is a complex number, then the real and imaginary parts
of $z$ are equal to $(z + \overline{z})/ 2$ and $(z - \overline{z})/(2
i)$, respectively.  The absolute values of the real and imaginary
parts of $z$ are less than or equal to the modulus of $z$.

	Let $z$, $w$ be complex numbers.  The modulus of the product
of $z$ and $w$ is equal to the product of their moduli, since the
complex conjugate of a product is equal to the product of complex
conjugates.  Furthermore,
\begin{eqnarray}
	|z + w|^2 & = & (z + w) (\overline{z} + \overline{w})	\\
		& = & |z|^2 + 2 \re z \, \overline{w} + |w|^2
							\nonumber \\
		& \le & |z|^2 + 2 |z| \, |w| + |w|^2 = (|z| + |w|)^2,
							\nonumber
\end{eqnarray}
so that $|z + w| \le |z| + |w|$.

	As in the case of real numbers, a sequence of complex numbers
$\{z_j\}_{j=1}^\infty$ converges to a complex number $z$ if for each
$\epsilon > 0$ there is an $L \ge 1$ such that
\begin{equation}
	|z - z_j| < \epsilon
\end{equation}
for all $j \ge L$.  This happens if and only if the sequences of real
and imaginary parts of the $z_j$'s converge to the real and imaginary
parts of $z$ as sequences of real numbers.

	Just as for real numbers, sums and products of convergent
sequences of complex numbers converge to the sum and product of the
corresponding limits, and the sequence of reciprocals of a convergent
sequence of nonzero complex numbers with nonzero limit converges to
the reciprocal of the limit.  The complex conjugates of a convergent
sequence of complex numbers converges to the complex conjugate of the
initial sequence.  One can verify these statements directly, or reduce
to the case of sequences of real numbers.

	If $\{z_j\}_{j=1}^\infty$ is a sequence of complex numbers
which converges to the complex number $z$, then
$\{|z_j|\}_{j=1}^\infty$ converges to $|z|$.  This follows from the
inequality
\begin{equation}
	\Bigl| |z| - |w| \Bigr| \le |z - w|
\end{equation}
for $z, w \in {\bf C}$, which can be derived from the triangle
inequality.

	A sequence of complex numbers $\{z_j\}_{j=1}^\infty$ is a
Cauchy sequence if for each $\epsilon > 0$ there is an $L \ge 1$ such
that
\begin{equation}
	|z_j - z_l| < \epsilon
\end{equation}
for all $j, l \ge L$.  It is easy to check that $\{z_j\}_{j=1}^\infty$
is a Cauchy sequence of complex numbers if and only if the real and
imaginary parts of the $z_j$'s form Cauchy sequences of complex
numbers.  As a result, every Cauchy sequence of complex numbers
converges.

	An infinite series $\sum_{j=1}^\infty a_j$ of complex numbers
converges if and only if the sequence of partial sums $\sum_{j=1}^n
a_j$ converges.  This happens if and only if the series of real and
imaginary parts of the $a_j$'s converge as series of real numbers.  We
say that $\sum_{j=1}^\infty a_j$ converges absolutely if
$\sum_{j=1}^\infty |a_j|$ converges.  This happens if and only if the
series of real and imaginary parts of the $a_j$'s converge absolutely.
If a series of complex numbers converges absolutely, then it
converges.

	If $p(z) = a_n \, z^n + \cdots + a_0$ is a polynomial function
on ${\bf C}$, with $n \ge 1$, $a_0, \ldots, a_n \in {\bf C}$, and $a_n
\ne 0$, then $p(z) = 0$ for at least one $z \in {\bf C}$.
Consequently such a polynomial can be expressed as the product of
$a_n$ times $n$ factors of the form $(z - z_j)$, with $z_1, \ldots,
z_n \in {\bf C}$.

\section{Quaternions}
\setcounter{equation}{0}

	By definition the quaternions form a $4$-dimensional
associative algebra over the real numbers denoted ${\bf H}$,
containing a copy of the real numbers which commute with all other
quaternions, and with the real number $1$ as the multiplicative
identity element for all of ${\bf H}$.  If $x \in {\bf H}$, then
\begin{equation}
	x = x_1 + x_2 \, i + x_3 \, j + x_4 \, k,
\end{equation}
for some $x_1, x_2, x_3, x_4 \in {\bf R}$, where $i, j, k \in {\bf H}$
satisfy
\begin{equation}
	i^2 = j^2 = k^2 = -1
\end{equation}
and
\begin{equation}
	i \, j = - j \, i = k,
\end{equation}
from which it follows that $i \, k = - k \, i = -j$ and $j \, k = k \,
j = i$.  If $x$ is as above, then we put
\begin{equation}
	x^* = x_1 - x_2 \, i - x_3 \, j - x_4 \, k,
\end{equation}
and define the modulus of $x$ to be the nonnegative real number
\begin{equation}
	|x| = \sqrt{x_1^2 + x_2^2 + x_3^2 + x_4^2}.
\end{equation}
One can check that $(x \, y)^* = y^* \, x^*$ for all $x, y \in H$ and
$|x|^2 = x \, x^* = x^* \, x$.  If $x$ is a nonzero quaternion, then
$x$ is invertible in ${\bf H}$ and $x^{-1} = x^* \, |x|^{-2}$.

	By the imaginary quaternions we mean the quaternions which are
real linear combinations of $i$, $j$, and $k$.  Every quaternion $x$
is the sum of a real number and an imaginary quaternion, which can be
expressed as $(x + x^*)/2$ and $(x - x^*)/2$.  If $w$ is an imaginary
quaternion, then $w^2 = - |w|^2$.  For all $x, y \in {\bf H}$ we have
that $|x \, y|^2 = (x \, y) (x \, y)^* = x \, y \, y^* \, x = x \,
|y|^2 \, x^* = |x|^2 \, |y|^2$ and therefore $|x \, y| = |x| \, |y|$.
One can check that $|x + y| \le |x| + |y|$ for all $x, y \in {\bf H}$.

\section{Real and complex vector spaces}
\setcounter{equation}{0}

	Let $V$ be a vector space over the real or complex numbers.
By a seminorm on $V$ we mean a nonnegative real-valued function
$N(v)$ defined for $v \in V$ such that $N(0) = 0$,
\begin{equation}
	N(\alpha \, v) = |\alpha| \, N(v)
\end{equation}
for all real or complex numbers $\alpha$, as appropriate, and all $v
\in V$, and
\begin{equation}
	N(v + w) \le N(v) + N(w)
\end{equation}
for all $v, w \in V$.  If $N(v) > 0$ for all $v \in V$ with $v \ne 0$,
then we say that $N$ is a norm on $V$.  Of course the usual absolute
value function or modulus are norms on ${\bf R}$, ${\bf C}$ as
one-dimensional vector spaces.

	Let $n$ be a positive integer, and consider the
$n$-dimensional vector spaces ${\bf R}^n$, ${\bf C}^n$ over ${\bf R}$,
${\bf C}$.  If $v$ is an element of ${\bf R}^n$ or ${\bf C}^n$ and $1
\le p < \infty$, put
\begin{equation}
	\|v\|_p = \bigg(\sum_{j=1}^n |v_j|^p \bigg)^{1/p},
\end{equation}
and when $p = \infty$ put
\begin{equation}
	\|v\|_\infty = \max \{|v_j| : 1 \le j \le n\}.
\end{equation}
It is easy to check directly that $\|v\|_p$ is a norm on ${\bf R}^n$,
${\bf C}^n$ when $p = 1, \infty$.  When $1 < p < \infty$, the
conditions except for the triangle inequality are very simple, and
that can be derived from the convexity of the function $t^p$ on the
nonnegative real numbers.  More precisely, one can use that to check
that the set of vectors $v$ with $\|v\|_p \le 1$ is convex, and then
derive the triangle inequality for $\|v\|_p$ from that.

	By inspection we have that
\begin{equation}
	\|v\|_\infty \le \|v\|_p
\end{equation}
for all $v$ in ${\bf R}^n$ or ${\bf C}^n$ and all $p$, $1 \le p <
\infty$.  Using this one can check that
\begin{equation}
	\|v\|_q \le \|v\|_p
\end{equation}
when $1 \le p \le q \le \infty$.  Similarly, it is easy to see that
\begin{equation}
	\|v\|_p \le n^{1/p} \, \|v\|_\infty
\end{equation}
when $1 \le p < \infty$.  In fact one has that
\begin{equation}
	\|v\|_p \le n^{(1/p) - (1/q)} \, \|v\|_q
\end{equation}
when $1 \le p \le q < \infty$.  This follows from the convexity of the
function $t^r$ on the nonnegative real numbers when $r \ge 1$.

	Suppose that $V$ is a real or complex vector space equipped
with a norm $\|v\|$.  A sequence $\{v_j\}_{j=1}^\infty$ of vectors in
$V$ is said to converge to $v \in V$ if for every $\epsilon > 0$
there is an $L \ge 1$ such that
\begin{equation}
	\|v - v_j\| < \epsilon
\end{equation}
for all $j \ge L$.  This is equivalent to $\lim_{j \to 0} \|v - v_j\|
= 0$ as a sequence of real numbers.  If $\{v_j\}_{j=1}^\infty$,
$\{w_j\}_{j=1}^\infty$ are sequences in $V$ which converge to $v, w
\in V$, respectively, then one can check that $\{v_j +
w_j\}_{j=1}^\infty$ converges to $v + w$.  If
$\{\alpha_j\}_{j=1}^\infty$ is a sequence of real or complex numbers
which converges to the real or complex number $\alpha$, as
appropriate, and it $\{v_j\}_{j=1}^\infty$ is a sequence of vectors in
$V$ which converges to $v \in V$, then $\{\alpha_j \,
v_j\}_{j=1}^\infty$ converges in $V$ to $\alpha \, v$.

	Using the triangle inequality one can check that
\begin{equation}
	\Bigl| \|v\| - \|w\| \Bigr| \le \|v - w\|
\end{equation}
for all $v, w \in V$.  It follows that if $\{v_j\}_{j=1}^\infty$ is a
sequence of vectors in $V$ which converges to a vector $v \in V$, then
$\{\|v_j\|\}_{j=1}^\infty$ converges to $\|v\|$ as a sequence of real
numbers.

	If $V$ is ${\bf R}^n$ or ${\bf C}^n$ with one of the norms
$\|\cdot \|_p$, $1 \le p \le \infty$, then convergence of vectors is
the same as for the standard Euclidean topologies.  In particular, a
sequence of vectors converges to another vector if and only if the $n$
sequences of coordinates converge to the corresponding coordinates of
the limit.  Actually this works for any norm on ${\bf R}^n$ or ${\bf
C}^n$.

	A finite-dimensional vector space $V$ over the real or complex
numbers is isomorphic to ${\bf R}^n$ or ${\bf C}^n$ for some $n$, and
there is a natural topology on $V$ so that any such isomorphism is a
homeomorphism.  One can show that any norm on $V$ induces the same
topology on $V$.

	A sequence $\{v_j\}_{j=1}^\infty$ of vectors in $V$ is said to
be a Cauchy sequence if for each $\epsilon > 0$ there is an $L \ge 1$
such that $\|v_j - v_l\| < \epsilon$ for all $j, l \ge L$.  Every
convergent sequence in $V$ is a Cauchy sequence.  If $V$ every Cauchy
sequence in $V$ converges, then $V$ is said to be a Banach space.

	If $V$ is ${\bf R}^n$ or ${\bf C}^n$ with any of the norms
$\|v\|_p$ described above, then a sequence in $V$ is a Cauchy sequence
if and only if the $n$ sequences of coordinates of the vectors are
Cauchy sequences, and it follows that every Cauchy sequence converges.
This holds for any finite-dimensional real or complex vector space
with respect to any norm, or by defining Cauchy sequences simply in
terms of the topology and the vector space structure.

	An infinite series $\sum_{j=1}^\infty v_j$ of vectors in a
real or complex vector space $V$ equipped with a norm $\|\cdot \|$ is
said to converge if the sequence of partial sums $\sum_{j=1}^n v_j$
converges.  An infinite series $\sum_{j=1}^\infty v_j$ of vectors in
$V$ converges absolutely if $\sum_{j=1}^\infty \|v_j\|$ converges as a
series of nonnegative real numbers.  If $\sum_{j=1}^\infty v_j$
converges absolutely, then one can check that the sequence of partial
sums is a Cauchy sequence.  If $V$ is complete, then the series
converges.  Conversely, if every absolutely convergent series in $V$
converges, then one can show that $V$ is complete.

	Suppose that $V_1$, $V_2$ are vector spaces, both real or both
complex, equipped with norms $\|\cdot \|_1$, $\|\cdot \|_2$,
respectively.  A linear mapping $T$ from $V_1$ to $V_2$ is said to be
bounded if there is a nonnegative real number $A$ such that
$\|T(v)\|_2 \le A \, \|v\|_1$ for all $v \in V_1$.  This is equivalent
to saying that $T$ is continuous as a mapping from $V_1$ to $V_2$ with
respect to the topologies defined by their norms.  Let us write
$\mathcal{BL}(V_1, V_2)$ for the space of bounded linear mappings from
$V_1$ to $V_2$.  Notice that this is a vector space over ${\bf R}$ or
${\bf C}$, as appropriate, because the sum of two bounded linear
mappings from $V_1$ to $V_2$ is a bounded linear mapping, and a scalar
multiple of a bounded linear mapping is bounded.

	If $V_1$ is ${\bf R}^n$ or ${\bf C}^n$ with one of the norms
described earlier, then every linear mapping from $V_1$ to $V_2$ is
bounded.  This is easy to see just using the definitions, writing a
vector in $V_1$ as a linear combination of the standard basis vectors.
More generally this works whenever $V_1$ has finite dimension.

	If $T$ is a bounded linear mapping from $V_1$ to $V_2$, then
its operator norm is denoted $\|T\|_{op, 12}$ and is defined to be the
supremum of $\|T(v)|_2$ over all $v \in V_1$ with $\|v\|_1 \le 1$.
Equivalently, the operator norm of $T$ satisfies the condition for the
nonnegative real number $A$ mentioned in the previous paragraph, and
is the smallest nonnegative real number with this property.  One can
check that the operator norm defines a norm on $\mathcal{BL}(V_1,
V_2)$.

	If $V_1$, $V_2$, $V_3$ are vector spaces, all real or all
complex, equipped with norms $\|\cdot \|_1$, $\|\cdot \|_2$, $\|\cdot
\|_3$, respectively, and if $T_1$, $T_2$ are bounded linear operators
from $V_1$ to $V_2$ and from $V_2$ to $V_3$, respectively, then the
composition $T_2 \circ T_1$ is a bounded linear operator from $V_1$ to
$V_3$.  The operator norm of $T_2 \circ T_1$ is less than or equal to
the product of the operator norms of $T_1$ and $T_2$.

	Suppose that $V_1$, $V_2$ are vector spaces, both real or both
complex, equipped with norms $\|\cdot \|_1$, $\|\cdot \|_2$,
respectively, and that $V_2$ is complete.  In this event one can show
that the vector space of bounded linear mappings from $V_1$ to $V_2$
is complete with respect to the operator norm.  If $W$ is a linear
subspace of $V_1$ which is dense in $V_1$, in the sense that every
element of $V_1$ can be expressed as the limit of a sequence of
vectors in $W$, and if $T$ is a linear mapping from $W$ to $V_2$ which
is bounded with respect to the restriction of the norm on $V_1$ to
$W$, then there is a unique extension of $T$ to a bounded linear
mapping from $V_1$ to $V_2$, with the same operator norm as the
original linear mapping on $W$.

\section{Real and complex algebras}
\setcounter{equation}{0}

	Let $\mathcal{A}$ be an algebra over the real or complex
numbers, and suppose that $\|\cdot \|$ is a norm on $\mathcal{A}$.

	Suppose that there is a nonnegative real number $C$ such that
\begin{equation}
	\|a \, b\| \le C \, \|a\| \, \|b\|
\end{equation}
for all $a, b \in \mathcal{A}$.  This implies that if
$\{a_j\}_{j=1}^\infty$, $\{b_j\}_{j=1}^\infty$ are sequences in
$\mathcal{A}$ which converge to $a, b \in \mathcal{A}$, then the
sequence of products $\{a_j \, b_j\}_{j=1}^\infty$ converges to the
product $a \, b$.  In fact this condition is equivalent to continuity
of the product at $0$, which is equivalent to continuity of the
product everywhere because of bilinearity.

	If $\mathcal{A}$ is finite-dimensional as a vector space
over the real or complex numbers, then this continuity condition is
automatic.

	Let $V$ be a nonzero real or complex vector space equipped
with a norm, and let $\mathcal{BL}(V)$ denote the algebra of bounded
linear operators on $V$, using composition as multiplication.  The
identity transformation $I$ on $V$ is the nonzero multiplicative
identity element, and its operator norm is equal to $1$.  If $T_1$,
$T_2$ are bounded linear operators on $V$, then the norm of the
composition $T_2 \circ T_1$ is less than or equal to the product of
the operator norms of $T_1$, $T_2$.  If $V$ is complete with respect
to its norm, then $\mathcal{BL}(V)$ is complete with respect to the
operator norm.

	Suppose now that $\mathcal{A}$ is an associative algebra over
the real or complex numbers with a nonzero multiplicative identity
element $e$.  Suppose further that $\|\cdot \|$ is a norm on
$\mathcal{A}$ such that
\begin{equation}
	\|e\| = 1
\end{equation}
and
\begin{equation}
	\|x \, y\| \le \|x\| \, \|y\|
\end{equation}
for all $x, y \in \mathcal{A}$, which means that $\mathcal{A}$ is a
normed algebra.  If one starts with a norm with the property that the
norm of a product is bounded by the product of the norms times a fixed
constant, then we can replace that norm with an equivalent one which
satisfies these conditions using the operator norm of the linear
transformation $x \mapsto a \, x$ on $\mathcal{A}$ with respect to the
initial norm on $\mathcal{A}$.  Let us ask also that $\mathcal{A}$ be
complete with respect to this norm, which is to say that $\mathcal{A}$
is a Banach algebra.

	Suppose that $x \in \mathcal{A}$ and that $\|x\| < 1$.  In
this event the series $\sum_{j=0}^\infty x^j$ converges absolutely in
$\mathcal{A}$ and hence converges.  As usual,
\begin{equation}
	(e - x) \, \bigg(\sum_{j=0}^n x^j \bigg)
		= \bigg(\sum_{j=0}^n x^j \biggr) (e - x)
		= e - x^{n+1}
\end{equation}
for all $n$.  It follows that $e - x$ is invertible in $\mathcal{A}$,
with inverse equal to $\sum_{j=0}^\infty x^j$.  More generally, for
each invertible element $w$ of $\mathcal{A}$ and each $x \in
\mathcal{A}$ with $\|x\| \, \|w^{-1}\| < 1$ we have that $w - x$ is an
invertible element of $\mathcal{A}$.

	In particular the set of invertible elements of $\mathcal{A}$
is an open subset of $\mathcal{A}$.  That is to say, for each
invertible element of $\mathcal{A}$ there is an open ball around that
element with respect to the norm on $\mathcal{A}$ which is contained
in the set of invertible elements of $\mathcal{A}$.  Using the series
expansion for $(e - x)^{-1}$ one can check that the mapping $w \mapsto
w^{-1}$ on the set of invertible elements of $\mathcal{A}$ is
continuous at $e$, and one can extend this to get continuity of the
multiplicative inverse at every invertible element of $\mathcal{A}$.

\section{Involutions}
\setcounter{equation}{0}

	Let $\mathcal{A}$ be an associative algebra over a field $k$.
An involution on $\mathcal{A}$ is a mapping $x \mapsto x^*$ which is
linear on $\mathcal{A}$ as a vector space over $k$ and which satisfies
\begin{equation}
	(x^*)^* = x
\end{equation}
for all $x \in \mathcal{A}$ and
\begin{equation}
	(x \, y)^* = y^* \, x^*
\end{equation}
for all $x, y \in \mathcal{A}$.  If $\mathcal{A}$ is an associative
algebra over a field $k$ with a nonzero multiplicative identity
element $e$ and an involution $x \mapsto x^*$, then $e^* = e$.  If $x$
is an invertible element of $\mathcal{A}$, then $x^*$ is an invertible
element of $\mathcal{A}$ too, with $(x^*)^{-1} = (x^{-1})^*$.  If
$\mathcal{A}$ is an algebra over the real or complex numbers equipped
with a norm, then it is natural to ask that an involution $x \mapsto
x^*$ be a bounded linear transformation with respect to the norm.
Frequently an involution will actually be an isometry.

	If $\mathcal{A}$ is commutative, then $x^* = x$ defines an
involution on $\mathcal{A}$.  This can be applied to the real numbers,
while complex conjugation defines an interesting involution on the
complex numbers as an algebra over the real numbers, and we saw
earlier that there is a natural involution on the quaternions as an
algebra over the real numbers.  Suppose that $\mathcal{A}$ is an
associative algebra over a field $k$, and let $M_n(\mathcal{A})$ be
the corresponding algebra of $n \times n$ matrices with entries in
$\mathcal{A}$.  If $x \mapsto x^*$ is an involution on $\mathcal{A}$,
then we get an involution on $M_n(\mathcal{A})$ by applying the
involution on $\mathcal{A}$ to each entry of the matrix and taking the
transpose of the matrix, i.e., interchanging the order of the indices.
If $\mathcal{A}$ is an associative algebra over a field $k$ with an
involution $x \mapsto x^*$, then the algebra $\mathcal{A}[t_1, \ldots,
t_n]$ of polynomials over the indeterminants $t_1, \ldots, t_n$
inherits a natural involution, by applying the involution on
$\mathcal{A}$ to the coefficients and defining $t_j^*$ to be $t_j$ for
$1 \le j \le n$.

	Let $\mathcal{A}$ be an associative algebra over a field $k$
equipped with an involution $x \mapsto x^*$.  For each $x, y \in
\mathcal{A}$ we have that
\begin{equation}
	(x \, y - y \, x)^* = - ( x^* \, y^* - y^* \, x^*).
\end{equation}
An element $w$ of $\mathcal{A}$ is said to be antisymmetric if $w^* =
- w$.  The antisymmetric elements of $\mathcal{A}$ form a linear
subspace of $\mathcal{A}$, which is closed under the operation $x \, y
- y \, x$.  Thus the antisymmetric elements of $\mathcal{A}$ form a
Lie subalgebra with respect to commutators, i.e., a Lie subalgebra of
$\lambda(\mathcal{A})$.

\section{Exponentiation}
\setcounter{equation}{0}

	The classical exponential mapping on the complex numbers ${\bf
C}$ can be defined by the power series
\begin{equation}
	\exp z = \sum_{n=0}^\infty \frac{z^n}{n!},
\end{equation}
where $n!$ is $n$ factorial, the product of the positive integers from
$1$ to $n$, which is interpreted as being equal to $1$ when $n = 0$.
When $n = 0$, $z^n$ is interpreted as being equal to $1$ for all
complex numbers $z$.

	Standard results in basic analysis such as the ratio test
imply that this series converges absolutely for all complex numbers
$z$.  If $z$ is a real number, then $z^n$ is a real number for all
nonnegative integers $n$, and $\exp z$ is a real number.  We also have
that
\begin{equation}
	\exp(z + w) = \exp(z) \, \exp(w)
\end{equation}
for all complex numbers $z$, $w$.

	Of course $\exp (0) = 1$, and hence
\begin{equation}
	\exp(z) \, \exp(-z) = 1.
\end{equation}
In particular, $\exp(z) \ne 0$ for all complex numbers $z$.

	If $x$ is a real number and $x \ge 0$, then $\exp x \ge 1$ by
inspection.  If $x \le 0$, then $0 < \exp x \le 1$ since $\exp(x) =
1/\exp (-x)$ and $\exp(-x) \ge 1$.  One can also check that the
exponential function is strictly increasing on the real line.

	If $z$ is a complex number then it follows easily from the
definition of the exponential function that
\begin{equation}
	\overline{\exp(z)} = \exp(\overline{z}).
\end{equation}
If $z = x + y \, i$, with $x, y \in {\bf R}$, then
\begin{equation}
	|\exp(z)|^2 = \exp (2 \, x).
\end{equation}

	By standard results in analysis the series expansion for the
exponential function converges uniformly on bounded subsets of ${\bf
R}$ or ${\bf C}$ and hence the exponential function is continuous.
The exponential function is actually differentiable of all orders.
The series expansion can be differentiated term by term, with the
well-known consequence that the derivative of the exponential function
is equal to itself.

	Now suppose that $\mathcal{A}$ is an associative algebra over
the real or complex numbers with a nonzero multiplicative identity
element $e$.  Suppose also that $\mathcal{A}$ is equipped with a norm
$\|\cdot \|$ which makes $\mathcal{A}$ a Banach algebra.  

	If $a \in \mathcal{A}$, then the exponential of $a$ is denoted
$\exp a$ and is defined in the same manner as before, as
$\sum_{n=0}^\infty a^n / n!$.  When $n = 0$ we interpret $a^n$ as
being equal to $e$ for all $a \in \mathcal{A}$.  The series converges
absolutely, since $\|a^n / n!\| \le \|a\|^n / n!$.  It follows in
particular that $\|\exp a\| \le \exp \|a\|$.  The series converges
uniformly on bounded subsets of $\mathcal{A}$ and therefore defines a
continuous mapping from $\mathcal{A}$ into itself.

	If $a, b \in \mathcal{A}$ commute, which is to say that $a \,
b = b \, a$, then $\exp (a + b)$ is equal to the product of $\exp a$
and $\exp b$, just as for real and complex numbers.  If we take $b = -
a$, then it follows that the product of $\exp a$ and $\exp (-a)$ is
equal to $e$.  Thus $\exp a$ is an invertible element of $\mathcal{A}$
for all $a \in \mathcal{A}$.

	Suppose that $\mathcal{A}$ is equipped with an involution $x
\mapsto x^*$ which is a bounded linear mapping on $\mathcal{A}$.  In
this event we have that $(\exp a)^*$ is equal to $\exp (a^*)$ for all
$a \in \mathcal{A}$.

	If $a \in \mathcal{A}$ is antisymmetric, $a^* = -a$, then
$(\exp a)^*$ is equal to $\exp (-a)$, the multiplicative inverse of
$\exp a$.  In any associative algebra over a field $k$ with an
involution $x \mapsto x^*$ and a mutliplicative identity element, the
invertible elements $x$ in the algebra with $x^{-1} = x^*$ form a
subgroup of the group of invertible elements with respect to
multiplication.  For a Banach algebra over the real or complex numbers
in which involution is a bounded linear mapping, the exponential
function sends the antisymmetric elements of the algebra into this
group.

\section{Power series}
\setcounter{equation}{0}

	Let $\mathcal{A}$ be an associative algebra over a field $k$,
and let $n$ be a positive integer.  We write $\mathcal{A}[[t_1,
\ldots, t_n]]$ for the algebra of formal power series in the
indeterminants $t_1, \ldots, t_n$ with coefficients in $\mathcal{A}$.
More precisely, each element of $\mathcal{A}[[t_1, \ldots, t_n]]$ can
be expressed as a formal sum
\begin{equation}
	\sum_\alpha a_\alpha \, t^\alpha,
\end{equation}
where the sum extends over all multi-indices $\alpha$, each $a_\alpha$
is an element of $\mathcal{A}$, and $t^\alpha$ is the monomial
associated to $\alpha$ discussed previously for polynomials.

	One can think of a power series in $n$ indeterminants $t_1,
\ldots, t_n$ and with coefficients in $\mathcal{A}$ as being defined
by a function from multi-indices into $\mathcal{A}$, which gives the
coefficients $a_\alpha$.  This makes precise the idea that the
monomials $t^\alpha$ are independent of each other in a simple way.

	One can view elements of $\mathcal{A}$ as power series in
which the coefficients of monomials of degree $\ge 1$ are equal to
$0$, and polynomials are the same as power series in which all but at
most finitely many terms are equal to $0$.  Therefore we have the
inclusions
\begin{equation}
	\mathcal{A} \subseteq \mathcal{A}[t_1, \ldots, t_n] 
		\subseteq \mathcal{A}[[t_1, \ldots, t_n]]
\end{equation}
in a natural way.

	One can add power series with coefficients in $\mathcal{A}$
and multiply them by elements of $k$ termwise, which makes
$\mathcal{A}[[t_1, \ldots, t_n]]$ a vector space over $k$.
Multiplication can be performed by grouping terms suitably, where the
coefficient of some monomial $t^\alpha$ in a product involves only
finitely many terms from the power series being multiplied, since
there are only finitely many pairs of multi-indices whose sum is equal
to $\alpha$.  This makes $\mathcal{A}[[t_1, \ldots, t_n]]$ an
associative algebra over $k$, where the indeterminants $t_1, \ldots,
t_n$ commute by definition, and which contains $\mathcal{A}$ and
$\mathcal{A}[t_1, \ldots, t_n]$ as subalgebras.  The linear operators
$\partial / \partial t_j$ can be defined on power series in the same
way as for polynomials, following the usual rules from calculus, and
are derivations on the algebra of power series.  If $\mathcal{A}$ is
commutative, then the algebra of power series with coefficients in
$\mathcal{A}$ is commutative too.

	Let us say that a sequence $\{p_j\}_{j=1}^\infty$ of power
series with coefficients in $\mathcal{A}$ in the indeterminants $t_1,
\ldots, t_n$ converges to a power series $p$ if for each multi-index
$\alpha$ there is a positive integer $L_\alpha$ such that the
coefficients of the monomial $t^\alpha$ in $p_j$ are equal to the
coefficient of $t^\alpha$ in $p$ when $j \ge L_\alpha$.  If
$\{p_j\}_{j=1}^\infty$, $\{q_j\}_{j=1}^\infty$ are sequences of power
series which converge to the power series $p$, $q$, then the sequences
$\{p_j + q_j \}_{j=1}^\infty$ and $\{p_j \, q_j \}_{j=1}^\infty$ of
sums and products of $p_j$'s and $q_j$'s converge to the sum $p + q$
and product $p \, q$ of the limits $p$, $q$, and if $c \in k$, then
$\{c \, p_j\}_{j=1}^\infty$ converges to $c \, p$.  Polynomials are
dense in the algebra of power series in the sense that every power
series can be expressed as a limit of a sequence of polynomials.  

	An infinite series $\sum_{j=1}^\infty r_j$ of power series
converges if the corresponding sequence of partial sums $\sum_{j=1}^n
r_j$ converges.  This happens if and only if the sequence of power
series $\{r_j\}_{j=1}^\infty$ converges to $0$, which is to say that
for each multi-index $\alpha$ the coefficient of $t^\alpha$ in $r_j$
is equal to $0$ for sufficiently large $j$.

	Let us suppose now that $\mathcal{A}$ contains a nonzero
multiplicative identity element $e$.  The power series with constant
term equal to $e$ and other coefficients equal to $0$ is the nonzero
multiplicative identity element in the algebra of power series with
coefficients in $\mathcal{A}$.

	Let $p$ be a power series with coefficients in $\mathcal{A}$
whose constant term is equal to $0$.  For each $j$ consider $p^j$, the
product $p \cdots p$ with a total of $j$ $p$'s, as a power series.
When $j = 0$ we can interpret $p^j$ as being the constant power series
$e$.  

	The coefficients of a monomial $t^\alpha$ in $p^j$ are equal
to $0$ when the degree of $\alpha$ is strictly less than $j$.  Hence
the sequence of $p^j$'s converges to $0$, and therefore
$\sum_{j=0}^\infty p^j$ converges.  By the usual computation, the
product of $e - p$ with $\sum_{j=0}^n p^j$ in either order is equal to
$e - p^{n+1}$, and it follows that $e - p$ has $\sum_{j=0}^\infty p^j$
as its multiplicative inverse.

	If a power series with coefficients in $\mathcal{A}$ has a
multiplicative inverse, then the constant term in the power series is
invertible as an element of $\mathcal{A}$.  Conversely, if a power
series with coefficients in $\mathcal{A}$ has invertible constant
term, then it is invertible as a power series.

	If $\mathcal{A}$ is equipped with an involution $a \mapsto
a^*$, then one can define $p^*$ for a power series $p$ by applying the
involution to the coefficients of $p$.  This defines an involution on
the algebra of power series with coefficients in $\mathcal{A}$.

\section{Exponentiation, 2}
\setcounter{equation}{0}

	Let $\mathcal{A}$ be an associative algebra over a field $k$
with characteristic $0$, and suppose that $\mathcal{A}$ has a nonzero
multiplicative identity element $e$.  Fix a positive integer $n$, and
let $p$ be a power series with coefficients in $\mathcal{A}$ in the
indeterminants $t_1, \ldots, t_n$ and with constant term equal to $0$.
Define the exponential of $p$ as a power series by
\begin{equation}
	\exp p = \sum_{n=0}^\infty \frac{1}{n!} \, p^n.
\end{equation}
Here the rational numbers $1/n!$ make sense as elements of $k$, and
hence as elements of $\mathcal{A}$ by taking multiples of $e$, because
$k$ is assumed to have characteristic $0$.  As usual we interpret
$p^n$ as being $e$ when $n = 0$.  Since the constant term of $p^n$ is
equal to $0$, the sequence of $p^n$'s converges to $0$ as a sequence
of power series.  Therefore the sum in $\exp p$ converges to a power
series with coefficients in $\mathcal{A}$.

	Notice that the constant term in $\exp p$ is equal to $e$ by
construction.  If $p$, $q$ are two power series with coefficients in
$\mathcal{A}$ whose constant terms are equal to $0$ and which commute
with each other, $p \, = q \, p$, then
\begin{equation}
	\exp (p + q) = (\exp p) (\exp q)
\end{equation}
by the standard computations.  In particular $\exp (-p)$ is the
multiplicative inverse of $\exp p$ in the algebra of power series with
coefficients in $\mathcal{A}$.

	If $\mathcal{A}$ is equipped with an involution $a \mapsto
a^*$, which induces an involution on the algebra of power series by
acting on the coefficients, and if $p$ is a power series with constant
term equal to $0$, then $p^*$ is a power series with constant term
equal to $0$ too, and
\begin{equation}
	\exp p^* = (\exp p)^*.
\end{equation}
If $p$ is antisymmetric in the sense that $p^* = -p$, then $h = \exp
p$ has the property that $h^* = h^{-1}$ in the algebra of power
series.

\section{$p$-Adic numbers}
\setcounter{equation}{0}

	Let $p$ be a prime number.  If $x$ is a rational number, then
the $p$-adic absolute value of $x$ is denoted $|x|_p$, and defined to
be $0$ when $x = 0$ and equal to $p^{-l}$ when $x = p^l \, a / b$,
where $a, b, l$ are integers, $a , b \ne 0$, and $a$, $b$ are not
integer multiples of $p$.

	It is easy to see that
\begin{equation}
\label{|x y|_p = |x|_p |y|_p}
	|x \, y|_p = |x|_p \, |y|_p
\end{equation}
for all $x, y \in {\bf Q}$.  One can also check that
\begin{equation}
\label{|x + y|_p le max(|x|_p, |y|_p)}
	|x + y|_p \le \max(|x|_p, |y|_p)
\end{equation}
for all $x, y \in {\bf Q}$.

	The $p$-adic numbers are denoted ${\bf Q}_p$ and are the
completion of the rational numbers with respect to the distance
function $|x - y|_p$.  More precisely, ${\bf Q}_p$ is a field which
contains a copy of ${\bf Q}$ as a subfield.  The $p$-adic absolute
value function $|x|_p$ is defined for all $x \in {\bf Q}_p$, with
$|x|_p = 0$ if and only if $x = 0$, and with (\ref{|x y|_p = |x|_p
|y|_p}), (\ref{|x + y|_p le max(|x|_p, |y|_p)}) valid for all $x, y
\in {\bf Q}_p$.

	A sequence $\{x_j\}_{j=1}^\infty$ of elements of ${\bf Q}_p$
converges to $x \in {\bf Q}_p$ if for each $\epsilon > 0$ there is an
$L \ge 1$ such that $|x - x_j|_p < \epsilon$ for all $j \ge L$, which
is the same as $\lim_{j \to \infty} |x - x_j|_p = 0$ as a limit of a
sequence of real numbers.  The rational numbers are dense in ${\bf
Q}_p$ in the sense that for each $x \in {\bf Q}_p$ there is a sequence
$\{x_j\}_{j=1}^\infty$ of rational numbers which converges to $x$ in
${\bf Q}_p$.

	A sequence $\{x_j\}_{j=1}^\infty$ in ${\bf Q}_p$ is a Cauchy
sequence if for each $\epsilon > 0$ there is an $L \ge 1$ such that
$|x_j - x_l|_p < \epsilon$ for all $j, l \ge L$.  Convergent sequences
are automatically Cauchy sequences, and ${\bf Q}_p$ is complete in the
sense that every Cauchy sequence in ${\bf Q}_p$ converges.

	If $x \in {\bf Q}_p$, then either $x = 0$ or $|x|_p$ is an
integer power of $p$, as one can show using the density of ${\bf Q}$
in ${\bf Q}_p$.  If $\{x_j\}_{j=1}^\infty$, $\{y_j\}_{j=1}^\infty$ are
sequences in ${\bf Q}_p$ which converge to $x, y \in {\bf Q}_p$, then
$\{x_j + y_j\}_{j=1}^\infty$ and $\{x_j \, y_j\}_{j=1}^\infty$
converge to $x + y$, $x \, y$, respectively.  If
$\{x_j\}_{j=1}^\infty$ is a sequence of nonzero elements of ${\bf
Q}_p$ which converges to $x \in {\bf Q}_p$, $x \ne 0$, then
$\{x_j^{-1}\}_{j=1}^\infty$ converges to $x^{-1}$.

	Because of the ultrametric version of the triangle inequality
(\ref{|x + y|_p le max(|x|_p, |y|_p)}), a sequence
$\{x_j\}_{j=1}^\infty$ in ${\bf Q}_p$ is a Cauchy sequence if $\lim_{j
\to \infty} x_j - x_{j+1} = 0$ in ${\bf Q}_p$.  Of course the converse
holds and works for real or complex numbers too.

	An infinite series $\sum_{j=1}^\infty a_j$ with terms in ${\bf
Q}_p$ converges if the sequence of partial sums $\sum_{j=1}^n a_j$
converges in ${\bf Q}_p$.  This happens if and only if $\lim_{j \to
\infty} a_j = 0$ in ${\bf Q}_p$.  For real or complex numbers
convergence of a series implies that the terms converge to $0$, and
the converse does not work in general.

	For ${\bf Q}_p$ one does not really need a separate notion of
absolute convergence.  For instance, if an infinite series
$\sum_{j=1}^\infty a_j$ converges in ${\bf Q}_p$, then every series
$\sum_{j=1}^\infty a_j \, b_j$ with $b_j \in {\bf Q}_p$, $|b_j|_p$
bounded, also converges in ${\bf Q}_p$, since $\lim_{j \to \infty} a_j
= 0$ implies that $\lim_{j \to \infty} a_j \, b_j = 0$ under these
conditions.  If $\sum_{j=1}^\infty a_j$ is an infinite series of real
or complex numbers such that $\sum_{j=1}^\infty a_j \, b_j$ converges
for all bounded sequences $\{b_j\}_{j=1}^\infty$ of real or complex
numbers, as appropriate, then $\sum_{j=1}^\infty a_j \, b_j$ converges
absolutely.

	If $x \in {\bf Z}$, then $|x|_p \le 1$.  One can check that if
$x \in {\bf Q}$ and $|x|_p \le 1$, then there is a sequence of
integers which converges to $x$ in ${\bf Q}_p$.

	Let ${\bf Z}_p$ denote the set of $x \in {\bf Q}_p$ such that
$|x|_p \le 1$.  This is the same as the set of $x \in {\bf Q}_p$ for
which there is a sequence of integers converging to $x$ in ${\bf
Q}_p$.  One can show that ${\bf Z}_p$ is a compact set in the sense
that every sequence of points in ${\bf Z}_p$ has a subsequence which
converges.

\section{${\bf Q}_p$ Vector spaces and algebras}
\setcounter{equation}{0}

	Let $p$ be a prime number, and let $V$ be a vector space over
${\bf Q}_p$.  A seminorm on $V$ is a function $N(v)$ defined for $v
\in V$ with values in the nonnegative real numbers such that
\begin{equation}
	N(\alpha \, v) = |\alpha|_p \, N(v)
\end{equation}
for all $\alpha \in {\bf Q}_p$ and $v \in V$ and
\begin{equation}
	N(v + w) \le N(v) + N(w)
\end{equation}
for all $v, w \in V$.  If
\begin{equation}
	N(v + w) \le \max(N(v), N(w))
\end{equation}
for all $v, w \in V$, then we say that $N$ is an ultrametric seminorm
on $V$.  A seminorm $N$ on $V$ is a norm if $N(v) = 0$ exactly when $v
= 0$, and an ultrametric norm is an ultrametric seminorm which is a
norm.

	Suppose that $N$ is a norm on $V$.  A sequence
$\{v_j\}_{j=1}^\infty$ of vectors in $V$ converges to a vector $v \in
V$ if for each $\epsilon > 0$ there is an $L \ge 1$ such that $N(v -
v_j) < \epsilon$ for all $j \ge L$.  Equivalently, $\lim_{j \to
\infty} v_j = v$ in $V$ if and only if $\lim_{j \to \infty} N(v - v_j)
= 0$ as a sequence of real numbers.  If $\{v_j\}_{j=1}^\infty$,
$\{w_j\}_{j=1}^\infty$ are sequences of vectors in $V$ which converge
to $v, w \in V$, then $\{v_j + w_j\}_{j=1}^\infty$ converges to $v +
w$.  If $\{\alpha_j\}_{j=1}^\infty$ is a sequence in ${\bf Q}_p$ which
converges to $\alpha \in {\bf Q}_p$ and $\{v_j\}_{j=1}^\infty$ is a
sequence in $V$ which converges to $v \in V$, then $\{\alpha_j \,
v_j\}_{j=1}^\infty$ converges to $\alpha \, v$ in $V$.

	For each positive integer $n$, the space ${\bf Q}_p^n$ of
$n$-tuples $x = (x_1, \ldots, x_n)$ with $x_j \in {\bf Q}_p$ for $1
\le j \le n$ is an $n$-dimensional vector space over ${\bf Q}_p$, with
respect to coordinatewise addition and scalar multiplication.  For $x
\in {\bf Q}_p^n$,
\begin{equation}
	\|x\| = \max (|x_1|_p, \ldots, |x_n|_p)
\end{equation}
defines an ultrametric norm on ${\bf Q}_p^n$.  A sequence in ${\bf
Q}_p^n$ converges to a vector in ${\bf Q}_p^n$ with respect to
$\|\cdot \|$ if and only if the coordinates of the vectors in the
sequence converge to the coordinates of the prospective limit in ${\bf
Q}_p$.

	Suppose now that $N$ is an ultrametric norm on a vector space
$V$ over ${\bf Q}_p$.  A sequence $\{v_j\}_{j=1}^\infty$ in $V$ is
said to be a Cauchy sequence if for each $\epsilon > 0$ there is an $L
\ge 1$ such that $N(v_j - v_l) < \epsilon$ for all $j, l \ge L$.
Because of the ultrametric version of the triangle inequality, this
happens if and only if $\lim_{j \to \infty} v_j - v_{j+1} = 0$ in $V$.
As usual every convergent sequence is a Cauchy sequence.  We say that
$V$ is complete if every Cauchy sequence in $V$ converges.

	An infinite series $\sum_{j=1}^\infty a_j$ with terms in $V$
converges if the sequence of partial sums $\sum_{j=1}^\infty a_j$
converges in $V$.  The sequence of partial sums is a Cauchy sequence
if and only if $\lim_{j \to \infty} a_j = 0$.  Every infinite series
$\sum_{j=1}^\infty a_j$ with terms in $V$ converges in $V$ if and only
if $V$ is complete.

	A sequence in ${\bf Q}_p^n$ is a Cauchy sequence with respect
to the norm mentioned previously if and only if its coordinates are
Cauchy sequences in ${\bf Q}_p$, and therefore every Cauchy sequence
in ${\bf Q}_p^n$ converges because of the completeness of ${\bf Q}_p$.

	Now suppose that $\mathcal{A}$ is an associative algebra over
${\bf Q}_p$ with a nonzero multiplicative identity element $e$ and a
norm $N$ such that $N(e) = 1$ and $N(a \, b) \le N(a) \, N(b)$ for all
$a, b \in \mathcal{A}$.  In particular if $\{a_j\}_{j=1}^\infty$,
$\{b_j\}_{j=1}^\infty$ are sequences in $\mathcal{A}$ which converge
to $a, b \in \mathcal{A}$, then $a_j \, b_j$ converges to $a \, b$ in
$\mathcal{A}$.

	Suppose also that $\mathcal{A}$ is complete.  If $x \in
\mathcal{A}$ and $N(x) < 1$, then $\lim_{n \to \infty} x^n = 0$ in
$\mathcal{A}$.  Hence $\sum_{j=0}^\infty x^j$ converges in
$\mathcal{A}$.  The product of the sum with $e - x$ in either order is
equal to $e$, which is to say that $e - x$ is invertible and $(e -
x)^{-1}$ is equal to $\sum_{j=0}^\infty x^j$.  More generally, if $x,
y \in \mathcal{A}$, $y$ is invertible in $\mathcal{A}$, and $N(x) \,
N(y^{-1}) < 1$, then $y - x$ is invertible in $\mathcal{A}$.

	As a basic example, consider the algebra of linear
transformations on ${\bf Q}_p^n$.  This can be identified with the
algebra of matrices $M_n({\bf Q}_p)$ with entries in ${\bf Q}_p$ in
the usual way.  Let us define the norm of a linear transformation on
${\bf Q}_p^n$ to be the maximum of $|a_{j, l}|_p$, $1 \le j, l \le n$,
where $(a_{j, l})$ is the corresponding matrix with entries in ${\bf
Q}_p$.  One can check directly that this norm has the properties
described in the previous paragraphs.  This norm is equal to the
operator norm of the linear transformation with respect to the norm
$\|x\| = \max (|x_1|_p, \ldots, |x_n|_p)$ on ${\bf Q}_p^n$.

	Now consider ${\bf Q}_p^n$, $n \ge 2$, with the norm which
assigns to $x = (x_1, \ldots, x_n)$ the maximum of $p^{-(j-1)/n} \,
|x_j|_p$ over $1 \le j \le n$.  Let $T$ be the linear mapping from
${\bf Q}_p^n$ to itself defined by $y = T(x)$ with $y_1 = p \, x_n$
and $y_j = x_{j-1}$, $1 \le j \le n-1$.  By construction, $T^n(x) = p
\, x$ for all $x \in {\bf Q}_p^n$.  The operator norm of $T$ with
respect to the norm just mentioned on ${\bf Q}_p^n$ is equal to
$p^{-1/n}$.

\section{${\bf Q}_p$ Exponentiation}
\setcounter{equation}{0}

	Fix a prime number $p$.  For each nonnegative integer $n$, we
would like to estimate $|n!|_p$.  Of course $n! = 1$ when $n = 0, 1$,
and hence $|n!|_p = 1$ too.  Basically we would like to estimate the
number of factors of $p$ in $n!$, $n \ge 2$.  The number of integers
from $1$ to $n$ which are divisible by $p$ is equal to $n / p$.  For
each positive integer $j$, the number of integers from $1$ to $n$
which are divisible by $p^j$ is equal to the integer part of $n /
p^j$.  The total number of factors of $p$ in $n!$ is equal to the sum
of these numbers.  In particular, it is strictly less than $\sum_{j =
1}^\infty n / p^j$.  Thus the total number of factors of $p$ in $n!$
is strictly less than $n / (p-1)$ when $n \ge 1$.  This implies that
$|1/n!|_p < p^{n/(p-1)}$ when $n \ge 1$.

	Suppose that $\mathcal{A}$ is an associative algebra over
${\bf Q}_p$ with a nonzero multiplicative identity element $e$ and an
ultrametric norm $N$ such that $N(e) = 1$, $N(a \, b) \le N(a) \,
N(b)$ for all $a, b \in \mathcal{A}$, and $\mathcal{A}$ is complete
with respect to $N$.  If $a \in \mathcal{A}$ and $N(a) <
p^{-1/(p-1)}$, then $N(a^n / n!) \le |1/n!|_p \, N(a)^n$ and $\lim_{n
\to \infty} a^n / n! = 0$ in $\mathcal{A}$.  We define $\exp a \in
\mathcal{A}$ to be the sum of the usual series $\sum_{n=0} a^n / n!$.

	Because of the ultrametric property for $N$, the set of $a \in
\mathcal{A}$ with $N(a) < p^{-1/(p-1)}$ is closed under addition.  If
$a, b \in \mathcal{A}$ have this property and $a$, $b$ commute in
$\mathcal{A}$, then $\exp (a + b)$ is equal to the product of $\exp a$
and $\exp b$.  By definition, $\exp 0 = e$.  If $a \in \mathcal{A}$
and $N(a) < p^{-1/(p-1)}$, then $\exp a$ is an invertible element of
$\mathcal{A}$, with inverse equal to $\exp (-a)$.  These remarks apply
in particular to $\mathcal{A} = {\bf Q}_p$ with $N(x) = |x|_p$.

\section{Traces and determinants}
\setcounter{equation}{0}

	Fix a positive integer $n$, and consider the algebra $M_n({\bf
C})$ of $n \times n$ matrices with complex entries.  This is an
associative algebra over the complex numbers with a nonzero
multiplicative identity element, given by the identity matrix.  We can
identify elements of $M_n({\bf C})$ with linear transformations on
${\bf C}^n$, and if we pick a norm on ${\bf C}^n$, such as the
standard Euclidean norm, then this leads to an operator norm on linear
transformations which makes $M_n({\bf C})$ into a Banach algebra.

	If $A$ is an $n \times n$ matrix with complex entries, the
exponential of $A$ is therefore defined.  A well-known theorem states
that the determinant of $\exp A$ is equal to the exponential of the
trace of $A$.

	This is trivial when $A$ is a diagonal matrix.  Similarly, it
is very easy to check this identity when $A$ is an upper-triangular
matrix.  One can derive the general case from this using the Jordan
canonical form.

	As another type of approach one can notice first that the
identity holds when $A$ is diagonalizable.  One can then argue that
this is a sufficiently large collection of matrices that the identity
should hold in general.

	For a third argument, consider $\exp (t \, A)$ as a
matrix-valued function on the real line, which is characterized by the
properties that it is equal to the identity matrix at $t = 0$ and
satisfies the differential equation that its derivative is equal to
$A$ times itself.  Of course the determinant of $\exp (t \, A)$ and
$\exp (t \tr A)$ are equal to $1$ at $t = 0$, and one can show that
they both satisfy the same differential equation.  This implies that
these functions are equal on the whole real line, and at $t = 1$ in
particular.

	Now suppose that $\mathcal{A}$ is an associative Banach
algebra over the real or complex numbers with a nonzero multiplicative
identity element $e$.  Let us consider the algebra $M_n(\mathcal{A})$
of $n \times n$ matrices with entries in $\mathcal{A}$.

	As a preliminary point let us consider norms on
$\mathcal{A}^n$ which make this into a Banach space and which are
compatible with multiplication and the norm on $\mathcal{A}$ in a nice
way.  For instance, using the norm on $\mathcal{A}$ we get a mapping
from $\mathcal{A}^n$ into the set of $n$-tuples of nonnegative real
numbers.  By applying one of the usual norms $\|\cdot \|_p$ on ${\bf
R}^n$, $1 \le p \le \infty$, we get a nice norm on $\mathcal{A}^n$.

	Each $n \times n$ matrix with entries in $\mathcal{A}$
corresponds to an operator on $\mathcal{A}^n$ in a natural way, and
with a suitable norm on $\mathcal{A}^n$ we get a nice operator norm on
$M_n(\mathcal{A})$.  Using such a norm $M_n(\mathcal{A})$ becomes a
Banach algebra, and it is nicer than that, with additional
compatibility with multiplication and the norm on $\mathcal{A}$.  At
any rate, we can consider exponentials on $M_n(\mathcal{A})$.  Let us
now suppose that $\mathcal{A}$ is commutative too.  Thus the
determinant and trace are defined as functions on $M_n(\mathcal{A})$
with values in $\mathcal{A}$, with their usual properties.

	Let us pause a moment and reflect on the determinant of the
exponential and the exponential of the trace.  For real or complex
matrices, these are scalar-valued analytic functions which can be
expressed by power series converging on the all of $M_n({\bf R})$ or
$M_n({\bf C})$.  The equality of these functions is equivalent to the
equality of the coefficients of the corresponding power series.  The
equality of the coefficients is a countable family of algebraic
identities, each involving only a finite number of terms with rational
coefficients.  The same identities imply that the determinant of the
exponential is equal to the exponential of the trace for $n \times n$
matrices with entries in a commutative Banach algebra.

	Let $k$ be a field, let $\mathcal{A}$ be an associative
algebra over $k$, and let $l$, $n$ be positive integers.
Notice that there is a natural equivalence
\begin{equation}
	M_n(\mathcal{A}[t_1, \ldots, t_l]) 
		\cong M_n(\mathcal{A})[t_1, \ldots, t_l].
\end{equation}
In other words, a matrix with entries in polynomials with coefficients
in $\mathcal{A}$ is basically the same as a polynomial with
coefficients in the algebra of matrices with entries in $\mathcal{A}$.
Similarly,
\begin{equation}
	M_n(\mathcal{A}[[t_1, \ldots, t_l]])
		\cong M_n(\mathcal{A})[[t_1, \ldots, t_l]],
\end{equation}
which is to say that matrices with entries in power series with
coefficients in $\mathcal{A}$ are basically the same as power series
with coefficients which are matrices with entries in $\mathcal{A}$.

	Now suppose that $k$ has characteristic $0$.  Suppose also
that $\mathcal{A}$ has a nonzero multiplicative identity element $e$,
which implies that $M_n(\mathcal{A})$ has a nonzero multiplicative
identity element given by the matrix with entries equal to $e$ on the
diagonal and to $0$ off of the diagonal.  Thus the exponential
function is defined for power series with constant term equal to $0$,
for power series with coefficients in $\mathcal{A}$ or
$M_n(\mathcal{A})$, where the result is a power series whose constant
term is the mutliplicative identity element.  Assuming that
$\mathcal{A}$ is commutative, we can again say that the determinant of
the exponential of a matrix whose entries are power series with
coefficients in $\mathcal{A}$ and constant term equal to $0$ is equal
to the exponential of the trace of the same matrix.  More precisely,
this statement follows from the same family of algebraic identities as
before.

	Let $p$ be a prime number.  Under suitable conditions we can
define exponentiation of elements of ${\bf Q}_p$ and of matrices with
entries in a commutative algebra over ${\bf Q}_p$, with restrictions
on the domain of the exponential in particular.  Formally the equality
between the determinant of the exponential of a matrix and the
exponential of the trace of the matrix follows from the usual family
of algebraic identities, with the extra ingredient now that the
relevant quantities be in the appropriate regions for the exponential
functions.

\end{document}